\documentclass[smallextended,natbib,authoryear
]{svjour3}       

\smartqed  

\usepackage{graphicx,color}

\usepackage{algorithm}
\usepackage{algorithmicx}
\usepackage{algpseudocode}

\usepackage{lineno}

\usepackage{mathptmx}      
\usepackage
{hyperref}

\journalname{Celestial Mechanics and Dynamical Astronomy}

\begin{document}

\title{ Intermediary LEO propagation including higher order zonal harmonics \thanks{Preliminary results were presented in KePASSA 2015, Toulouse, France, Oct.~28--30, 2015}
 }
\author{Denis Hautesserres \&
Martin Lara
}


\institute{D.~Hautesserres \at
            Centres de Comp\'etence Technique, Centre National d'\'Etudes Spatiales, \\
              Tel.: +34-856-119449\\
              \email{denis.hautesserres@cnes.fr}
\and
M.~Lara \at
           GRUCACI, University of La Rioja, and Space Dynamics Group -- UPM, \\
              Tel.: +34-856-119449\\
              \email{mlara0@gmail.com}
}

\date{\color{red}DRAFT PAPER of September 19, 2016}

\maketitle

\begin{abstract}
{Two} new intermediary orbits of the artificial satellite problem are proposed. The analytical solutions include higher order effects of the Geopotential, and are obtained by means of a torsion transformation applied to the quasi-Keplerian system resulting after the elimination of the parallax simplification, {for the first intermediary, and after the elimination of the parallax and perigee simplifications, for the second one}. The new intermediaries perform notably well for low earth orbits propagation, are free from special functions, and result advantageous, both in accuracy and efficiency, when compared to the standard Cowell integration of the $J_2$ problem, thus providing appealing alternatives for onboard, short-term, orbit propagation under limited computational resources.
\end{abstract}




\section{Introduction}

Intermediary orbits provide approximate analytical solutions to the artificial satellite problem. By definition, \emph{common}  intermediaries of the main problem, which are simplifications of the Geopotential where only the $J_2$ term is taken into account, must capture all the secular effects up to the first order of $J_2$ \citep{GarfinkelAksnes1970}. In the \emph{natural} perspective, main problem intermediaries are able to cope with all the first order periodic effects in addition to the secular terms \citep{Deprit1981}.
\par

In spite of their limited precision, it has recently been demonstrated that evaluation of intermediary orbits can compete in some cases with standard numerical integration of simplified models, and, in particular, in the case of short-term onboard orbit propagation of low earth orbits (LEO). Indeed, the intrinsic uncertainty with which the satellite position is usually known onboard makes the errors obtained with both approaches, numerical integration and intermediary evaluation, to share similar statistics \citep{GurfilLara2014CeMDA}. This characteristic motivates the current interest in the search for new approximate analytical solutions of the artificial satellite problem \citep{MartinusiDell-ElceKerschen2015}. Besides, the fact that, contrary to usual intermediaries, Deprit's (\citeyear{Deprit1981}) radial intermediary  does not depend on special functions, makes its evaluation straightforward and, therefore, proves to be a suitable candidate for the implementation of orbit propagators to be used onboard. 
\par

On the other hand, {exclusion of higher order harmonics is known to bring relevant errors into the predicted orbit. Besides},
it is known that the uncertainties introduced in the computation of the initial conditions in the intermediary variables when neglecting periodic effects of the second order of $J_2$, are an important source of errors in the short-term propagation of artificial satellites \citep{Lara2015ASR}. Hence, inclusion of as many higher order short-period effects as possible, in addition to the secular terms of the second order, seems imperative for extending the validity of intermediary solutions to encompass several orbital periods or days.
\par

We explore the performances of natural intermediary solutions that include the gravitational effects of the earth's zonal harmonics up to $J_4$. Like in the case of the main problem, the intermediaries are obtained from the simplified Hamiltonian stemming from the elimination of the parallax \citep{Deprit1981,LaraSanJuanLopezOchoa2013b}. {In particular, two radial intermediaries are discussed which show similar efficiency for short-term propagation, yet the second intermediary improves the accuracy by accounting for long-period terms.}
\par

The first intermediary is obtained by neglecting the dynamics associated to the trigonometric terms in the argument of the perigee; the resulting Hamiltonian is radial and, therefore, can be integrated by separation of variables. Alternatively, it is converted into a pure Keplerian system by means of a torsion transformation \citep{DepritJGCD1981}. Because the neglected, second order terms are factorized by the eccentricity, the use of the intermediary is constrained to the case of the lower eccentricity orbits, which is the case of most operational satellites in LEO. This intermediary includes the secular terms of the solution up to the second order of $J_2$ and the short-period corrections up to the same order. Due to the precision with which the short-period corrections are computed, the uncertainties in the computation of initial conditions in the intermediary variables are very small, in this way enlarging the time span in which ephemeris predictions are of acceptable accuracy.
\par

However, including second order short-period corrections in the intermediary is at the cost of degrading computation performance in terms of time and power consumption. Since the intermediary is intended for onboard orbit propagation, where both memory allocation and power consumption may be an issue, several simplifications to reduce the computational burden to a minimum are applied. In particular, a simplified version of the intermediary is developed in which the direct short-period corrections are limited to first order effects, while keeping the secular terms up to the second order. Still, the inverse short-period corrections include the most relevant second order terms, and hence the simplified intermediary has comparable precision to the full version in the computation of initial conditions in the intermediary variables chart. This simplified version of the intermediary clearly provides better accuracy than other main problem intermediaries, and in particular than Deprit's (first order) radial intermediary \citep{GurfilLara2014CeMDA}, whereas the computational effort is only increased slightly. Therefore, the simplified version of the intermediary can advantageously compete with the usual Cowell integration of the $J_2$ problem for short-term prediction.
\par

{The second intermediary is designed to improve accuracy. Indeed,} the first radial intermediary has been constructed by ignoring those second order terms of the Hamiltonian after the elimination of the parallax which are factored by the eccentricity. As a consequence, the intermediary solution is affected by long-period errors related to the neglected perigee dynamics. Even though the amplitude of these long-period errors is small and the missing terms do not deteriorate too much the intermediary solution in the short-time spans for which the intermediary predictions are intended, their effects are clearly apparent since the beginning of the propagation. {On the contrary, the second intermediary recovers} the neglected long-period effects by carrying out an additional transformation. In the wake of classic analytical perturbation methods, the long-period dynamics can be removed, rather than ignored, by means of the elimination of the perigee simplification \citep{AlfriendCoffey1984,LaraSanJuanLopezOchoa2013c}. Regrettably, the additional equations introduced in the solution by this new canonical transformation slow down the evaluation of the second intermediary, on one side, and give rise to the small divisors problem related to the critical inclination resonance, on the other. However, both inconveniences are easily avoided in the case of LEO, where second order terms factored by the \emph{square} of the eccentricity, which appear associated to the even zonal harmonics, can be directly neglected. The remaining long-period terms, which are related to the latitudinal asymmetry of the Geopotential and do not experience the small divisors problem, are recovered by means of extremely simple corrections of straightforward evaluation. In addition, the transformation equations of the elimination of the perigee of the terms related to the $J_3$ harmonic coefficient can be dramatically simplified in the case of LEO, in this way producing effective corrections while insignificantly increasing the evaluation effort of the second intermediary solution. 
\par

It is not a surprise that the simplified long-period gravitational corrections used by the second intermediary match corresponding ones in the SGP4 propagator \citep{HootsRoehrich1980,ValladoCrawfordHujsakKelso2006}, even though the latter come from Brouwer's gravitational solution in which the perigee is averaged by a canonical transformation only \emph{after} the short-period terms have been completely removed from the original Hamiltonian \citep{Brouwer1959}. Quite on the contrary, the essential short-period terms of the zonal problem remain in the quasi-Keplerian, intermediary solution. However, since the elimination of the perigee is carried out only up to the first order, this identity was indeed expected from the commutativity of infinitesimal transformations.
\par

The paper is organized as follows. First, the dynamical model used is described in Section \ref{s:dm}. It follows in Section \ref{s:lis} the derivation of the first intermediary solution using the usual approach. Then, Section \ref{s:ne} presents some numerical tests that show the {higher accuracy} of the new intermediary over the numerical integration of the main problem. Next, allowable simplifications that speed up evaluation of the first intermediary while keeping sufficient accuracy are discussed in Section \ref{s:fe}. The improvements achieved by the second intermediary, where long-period corrections are taken into account, are analyzed in Section \ref{s:lpc}. { Finally, the procedure for obtaining the analytical solutions as well as the algorithms that implement them are recapitulated  in Section \ref{s:recapitulation}, where runtime comparisons between the $J_2$ numerical integration and the analytical solutions are also carried out for the particular case of Planet Labs-Dove satellites.\footnote{\href{http://www.nasa.gov/mission_pages/station/research/experiments/1326.html}{www.nasa.gov/mission\_pages/station/research/experiments/1326.html}. Accessed: Sep\-tem\-ber 10, 2016.}}
\par

The incorporation of a realistic atmospheric drag model into the intermediary solution, whose effect can be as important as the $J_2$ disturbances for the lower altitude orbits, is not discussed in the paper, which constrains its scope to providing efficient alternatives for those problems where the Cowell integration of the $J_2$ problem is the standard approach.

\section{Dynamical model} \label{s:dm}

The model includes effects of the main gravitational harmonics of the Geopotential, namely those corresponding to the zonal harmonic coefficients $J_2$, $J_3$, and $J_4$. Missing other disturbing effects, as atmospheric drag, lunisolar perturbations or tesseral and higher order zonal harmonics, of course makes predictions based on the $J_2$--$J_4$ model of limited precision. However, this simple model may be useful for short-term propagation, say in time spans ranging from a few minutes to a few orbits, where it would improve predictions in all those cases in which the simpler $J_2$ model is currently been used.
\par

Since the dynamical model is conservative the problem is formulated in the Hamiltonian frame,
\begin{equation} \label{zonal}
\mathcal{H}=\mathcal{H}_\mathrm{K}+\mathcal{Z},
\end{equation}
where $\mathcal{H}_\mathrm{K}$ represents the Keplerian attraction, and $\mathcal{Z}$ is the disturbing function which encompasses the non-centralities of the gravitational potential due to the zonal harmonics { \citep[see, for instance,][ch.~5]{Kellogg1953}}, viz.
\begin{equation} \label{distpot}
\mathcal{Z}=\frac{\mu}{r}\sum_{m\ge{2}}\frac{\alpha^m}{r^m}C_{m,0}\,P_m(\sin\varphi),
\end{equation}
where $\mu$ is the earth's gravitational parameter, $r$ is the distance from the earth's center of mass, $\varphi$ is latitude, the scaling factor $\alpha$ is the earth's equatorial radius, $P_m$ are Legendre polynomials of degree $m$, and $C_{m,0}=-J_m$ are corresponding zonal harmonic coefficients. For the earth, $C_{2,0}$ is of the order of one thousandth, while zonal coefficients of higher degree than the second are of the order of one millionth.
\par

The Keplerian part of the Hamiltonian is written as
\begin{equation} \label{KeplerianDelaunay}
\mathcal{H}_\mathrm{K}=-\frac{\mu}{2a},
\end{equation}
where $a$ is the orbit semi-major axis. Besides, in view of $\sin\varphi=\sin{I}\sin\theta$, where $I$ is orbital inclination and $\theta$ the argument of the latitude, up to degree 4, the Legendre polynomials in Eq.~(\ref{distpot}) are written
\begin{eqnarray*}
P_2 &=& -\frac{1}{2}+\frac{3}{4}s^2-\frac{3}{4}s^2\cos2\theta, \\
P_3 &=& \left(-\frac{3}{2}+\frac{15}{8}s^2\right)s\sin\theta-\frac{5}{8}s^3\sin3\theta, \\
P_4 &=& \frac{3}{64} \left(8-40s^2+35s^4\right)+\frac{5}{16} \left(6-7s^2\right)s^2\cos2\theta+\frac{35}{64} s^4 \cos4\theta,
\end{eqnarray*}
where the abbreviation $s$ stands for the sine of the inclination. Furthermore, since we are using Hamiltonian formalism, we assume that $a$, $r$, $\theta$, and $s$ are functions of the Delaunay canonical variables, which are given by the mean anomaly $\ell$, the argument of perigee $g$, the right ascension of the ascending node $h$, the Delaunay action $L$, the modulus of the angular momentum $G$, and the projection of the angular momentum vector along the earth's rotation axis $H$. In particular,
\begin{equation} \label{radius}
a=\frac{L^2}{\mu}, \qquad
s=\sqrt{1-c^2}, \quad c=\frac{H}{G}, \qquad
\theta=f+g, \qquad
r=\frac{p}{1+e\cos{f}},
\end{equation}
where the conic parameter $p$ and the orbit eccentricity $e$ are
\begin{equation}
p=\frac{G^2}{\mu}, \qquad
e=\sqrt{1-\eta^2}, \quad \eta=\frac{G}{L}, 
\end{equation}
and $f\equiv{f}(\ell,L,G)$, the true anomaly, is an implicit function of $\ell$ whose computation requires solving the Kepler equation.
\par

The zonal Hamiltonian given by Eqs.~(\ref{zonal}) and (\ref{distpot}) is a two degrees of freedom (DOF) Hamiltonian because, due to the axial symmetry of the zonal problem, the right ascension of the ascending node is a cyclic variable; still, the general solution to the flow derived from Eqs.~(\ref{zonal}) and (\ref{distpot}) is not known. However, in some cases, the zonal Hamiltonian can be simplified to the extent of providing useful approximate solutions, whose validity is constrained to certain regions of phase space and apply only during a limited time interval. This is the case of the LEO region, where the eccentricity is commonly small and, as a consequence, some of the effects related to the perigee dynamics only accumulate to the extent of being observable after a long time.
\par

\section{LEO Intermediary solution} \label{s:lis}

The original Hamiltonian in Eq.~(\ref{zonal}) is first simplified by carrying out the elimination of the parallax simplification, which is a canonical transformation $(\ell,g,h,L,G,H)\rightarrow(\ell',g',h',L',G',H')$, from old (osculating) to new (prime) variables, that removes non-essential short-period effects from the Hamiltonian \citep{Deprit1981}. {Second order terms of the $J_2$ Hamiltonian after the elimination of the parallax have been repeatedly published in the literature, including Deprit's original paper. In the case of higher order zonal harmonics, after reformulating corresponding terms of Eq.~\ref{distpot} using the identity $r^{-n}=r^{-2}p^{n-2}(1+e\cos{f})^{2-n}$, they are trivially obtained up to the second order of $J_2$ by dropping terms that explicitly depend on the true anomaly \citep[see][and references therein for further details]{LaraSanJuanLopezOchoa2013c,LaraSanJuanLopezOchoa2013b}}.
\par

After truncation to the second order of $J_{2}$, we get the Hamiltonian in new variables $\mathcal{H}_1\equiv\mathcal{H}_1(\ell',g',-,L',G',H')$
\begin{eqnarray} \label{Kfull}
\mathcal{H}_1 &=& -\frac{\mu}{2a} -\frac{\mu}{2p}\frac{\alpha^2}{r^2}J_{2}\left(1-\frac{3}{2} s^2\right)
+\frac{\mu}{2p}\frac{p^2}{r^2}J_{2}^2\Bigg\{\frac{1}{4}\frac{\alpha^4}{p^4}\Bigg[\frac{1}{4}-\frac{21}{4}c^4  \\ \nonumber
&& +3 \tilde{J}_4 \left(1-5 s^2+\frac{35}{8}s^4\right)\Bigg]
+\frac{3}{4}\frac{\alpha^3}{p^3}\tilde{J}_3\left(1-5c^2\right)es\sin\omega  \Bigg\}
+\mathcal{O}(e^2J_{2}^2),
\end{eqnarray}
where
\begin{equation}
\tilde{J}_3=\frac{J_{3}}{J_{2}^2}, \qquad \tilde{J}_4=\frac{J_{4}}{J_{2}^2},
\end{equation}
are of order 1, and now $a$, $p$, $e$, $r$, $s$, and $c$, are functions of the prime variables, while $\omega=g'$.
\par

Neglecting $\mathcal{O}(e^2J_{2}^2)$ is a reasonable assumption in the case of common operational satellites in LEO because of their low eccentricities. Then, effects of the perigee dynamics remain limited to the contribution of the zonal harmonic of degree 3, which is $\mathcal{O}(eJ_{2}^2)$. If we further ignore this term in Eq.~(\ref{Kfull}), the argument of the perigee no longer appears in the Hamiltonian, which, by this reason, turns out to be integrable. We remark that integrability is only reached at the cost of missing essential effects related to the long-period dynamics, which, therefore, constrains applicability of the solution to some regions of phase space. In particular, the removal of the perigee either by truncation or by perturbation theory prevents application of the solution to the case of inclination resonances \citep[see][and references therein]{Lara2015IR}.
\par

The integrability of the truncated Hamiltonian up to $\mathcal{O}(eJ_{2}^2)$ is better shown in the chart of polar-nodal, canonical variables $(r,\theta,\nu,R,\Theta,N)$, standing from radial distance, argument of the latitude, right ascension of the ascending node, radial velocity, total angular momentum, and polar component of the angular momentum along the earth's rotation axis, respectively. The transformation from Delaunay to polar-nodal variables is given by
\begin{equation}
r=\frac{p}{1+e\cos{f}}, \quad
\theta=f+g, \quad
\nu=h, \quad
R=\frac{G}{p}\,e\sin{f}, \quad
\Theta=G, \quad 
N=H,
\end{equation}
where the two first equations were previously given in Eq.~(\ref{radius}).
\par

Then, after neglecting in Eq.~(\ref{Kfull}) those terms which are factorized by the eccentricity, the simplified Hamiltonian is written
\begin{eqnarray} \label{DRI}
\mathcal{K} &=& \frac{1}{2}\left(R^2+\frac{\Theta^2}{r^2}\right)-\frac{\mu}{r}
+\frac{1}{2}\frac{\Theta^2}{r^2}\,\epsilon\left(2-3s^2\right)
\\ \nonumber
&& 
+\frac{1}{2}\frac{\Theta^2}{r^2}\,\epsilon^2\Bigg[\left(\frac{1}{4}-\frac{21}{4} c^4\right)
+ 3\tilde{J}_{4}\left(1-5 s^2+\frac{35}{8}s^4\right)
\Bigg],
\end{eqnarray}
where we abbreviate
\begin{equation} \label{epsilon}
\epsilon=-\frac{1}{2}\frac{\alpha^2}{p^2}J_{2},
\end{equation}
and now, $c$ and $p$ are functions of the polar-nodal variables, to wit $c=N/\Theta$ and $p=\Theta^2/\mu$. 
Note that, because the argument of the perigee has been ignored in Eq.~(\ref{Kfull}), $\theta$ is cyclic in Eq.~(\ref{DRI}) and, in consequence, $\Theta$ is constant and so it is $\epsilon\equiv\epsilon(\Theta)$. 
\par

Equation (\ref{DRI}) can be reorganized as
\begin{equation} \label{quasiK}
\mathcal{K} = \frac{1}{2}\left(R^2+\frac{\Theta^2}{r^2}\Phi^2\right)-\frac{\mu}{r},
\end{equation}
where
\begin{equation} \label{Phi2}
\Phi^2\equiv\Phi(\Theta,N)^2=1-\epsilon\left(1-3c^2\right)+\frac{1}{4}\epsilon^2\left[1-21c^4
+ \frac{3}{2}\tilde{J}_{4}\left(3-30c^2+35c^4\right)
\right],
\end{equation}
is also constant. The new Hamiltonian (\ref{quasiK}) represents a \emph{quasi-Keplerian} system with modified ``angular momentum''
$\tilde\Theta=\Theta\,\Phi(\Theta,N)$. It can be further reduced to a pure Keplerian system by means of a \emph{torsion} transformation $(r,\theta,\nu,R,\Theta,N)\longrightarrow(\tilde{r},\tilde\theta,\tilde\nu,\tilde{R},\tilde\Theta,\tilde{N})$,
which is a canonical transformation that modifies the angular variables while leaving untouched the radial ones \citep{Deprit1981,DepritJGCD1981}. 
\par

In the style of the Hamilton-Jacobi equation method, the torsion is defined by the transformation in mixed variables
\begin{equation}
r=\frac{\partial{T}}{\partial{R}}, \quad \theta=\frac{\partial{T}}{\partial\Theta}, \quad \nu=\frac{\partial{T}}{\partial{N}},
\qquad
\tilde{R}=\frac{\partial{T}}{\partial\tilde{r}}, \quad \tilde\Theta=\frac{\partial{T}}{\partial\tilde\theta}, \quad \tilde{N}=\frac{\partial{T}}{\partial\tilde\nu},
\end{equation}
which is derived from the generator
\begin{equation}
T\equiv{T}(\tilde{r},\tilde\theta,\tilde\nu,{R},\Theta,{N})=\tilde{r}R+\tilde\theta\,\Theta\,\Phi(\Theta,N)+\tilde\nu{N}.
\end{equation}
Then, the transformation equations are $\tilde{r}=r$, $\tilde{R}=R$, $\tilde{N}=N$, and
\begin{equation} \label{transform0}
\theta=\tilde\theta\left(\Phi+\Theta\frac{\partial\Phi}{\partial\Theta}\right), 
\qquad 
\nu=\tilde\nu+\tilde\theta\,\Theta\frac{\partial\Phi}{\partial{N}},
\qquad 
\tilde\Theta=\Theta\,\Phi.
\end{equation}
\par

{Using the chain rule, Eq.~(\ref{transform0}) is rewritten as
\begin{equation} \label{transform}
\theta=\frac{\tilde\theta}{\Phi}\left(\Phi^2-2\epsilon\frac{\partial\Phi^2}{\partial\epsilon}-\frac{1}{2}c\,\frac{\partial\Phi^2}{\partial{c}}\right),
\qquad 
\nu=\tilde\nu+\frac{1}{2}\frac{\tilde\theta}{\Phi}\,\frac{\partial\Phi^2}{\partial{c}},
\qquad 
\tilde\Theta=\Theta\,\Phi,
\end{equation}
where
\begin{eqnarray} \label{dPhidc}
\frac{\partial\Phi^2}{\partial{c}} &=& 
3\epsilon\,c  \left\{2-\epsilon \left[7 c^2+\frac{5}{2}\left(3-7 c^2\right)\tilde{J}_{4}\right]\right\}, \\  \label{dPhidep}
\frac{\partial\Phi^2}{\partial\epsilon} &=&
-1+3 c^2+\frac{1}{2}\epsilon\left[1-21c^4+\frac{3}{2}\left(3-30c^2+35c^4\right)\tilde{J}_{4}\right].
\end{eqnarray}
}
\par

Note that, while the transformation from old to new (tilde) variables is completely explicit, computing the transformation from new to old variables requires solving $\Theta$ from the implicit equation $\Theta\,\Phi(\Theta,\tilde{N})-\tilde\Theta=0$, as derived from the last of Eq.~(\ref{transform}), {in which it is simple to check from Eqs.~(\ref{Phi2}) and (\ref{epsilon}) that $\Phi^2$ can be written as a polynomial of degree 6 in $c^2$}. A single Newton-Raphson iteration of this equation is enough to provide the required accuracy. Indeed, neglecting terms of the order of $J_{2}^3$ and higher, we get
\begin{equation} \label{ThetaNR}
\Theta=\tilde\Theta\left\{1+\frac{1}{2}\tilde\epsilon\left(1-3\tilde{c}^2\right)-
\frac{3}{4}\tilde\epsilon^2\left[\frac{1}{4}\left(3-30\tilde{c}^2+35\tilde{c}^4\right)\tilde{J}_{4}+1-7\tilde{c}^2+10\tilde{c}^4\right]\right\},
\end{equation}
where
\begin{equation}
\tilde{c}=\frac{\tilde{N}}{\tilde\Theta}, \qquad \tilde\epsilon=-\frac{1}{2}\frac{\alpha^2}{\tilde{p}^2}J_{2}, \qquad \tilde{p}=\frac{\tilde\Theta^2}{\mu}.
\end{equation}
\par

Finally, the quasi-Keplerian Hamiltonian (\ref{quasiK}) is written in the new variables, to give
\begin{equation} \label{tildeKepler}
\mathcal{Q}\equiv\mathcal{K}(\tilde{r},\tilde\theta,\tilde\nu,\tilde{R},\tilde\Theta,\tilde{N})
=\frac{1}{2}\left(\tilde{R}^2+\frac{\tilde\Theta^2}{\tilde{r}^2}\right)-\frac{\mu}{\tilde{r}},
\end{equation}
which is a true Keplerian system in the chart $(\tilde{r},\tilde\theta,\tilde\nu,\tilde{R},\tilde\Theta,\tilde{N})$, and whose solution is standard.
\par

In summary, a compact analytical solution of the zonal problem has been computed, which neglects secular and periodic effects of $\mathcal{O}(e^2J_2^2)$ as well as long-period effects of $\mathcal{O}(eJ_2^2)$. Therefore, this intermediary solution should be accurate enough for short-term propagation of the lower eccentricity orbits in the LEO region.

\section{Numerical tests} \label{s:ne}

In order to check the usefulness of the intermediary solution, for a variety of test cases we compare it with the numerical, {Runge-Kutta (R-K)} integration of the original problem in Cartesian coordinates ---the flow derived from Eqs.~(\ref{zonal}) and (\ref{distpot}) with $m=4$, $\sin\varphi=z/r$, and $r=\sqrt{x^2+y^2+z^2}$, hereafter called the full zonal model.
\par

\begin{figure}[htbp]
\centerline{
\includegraphics[scale=0.8]{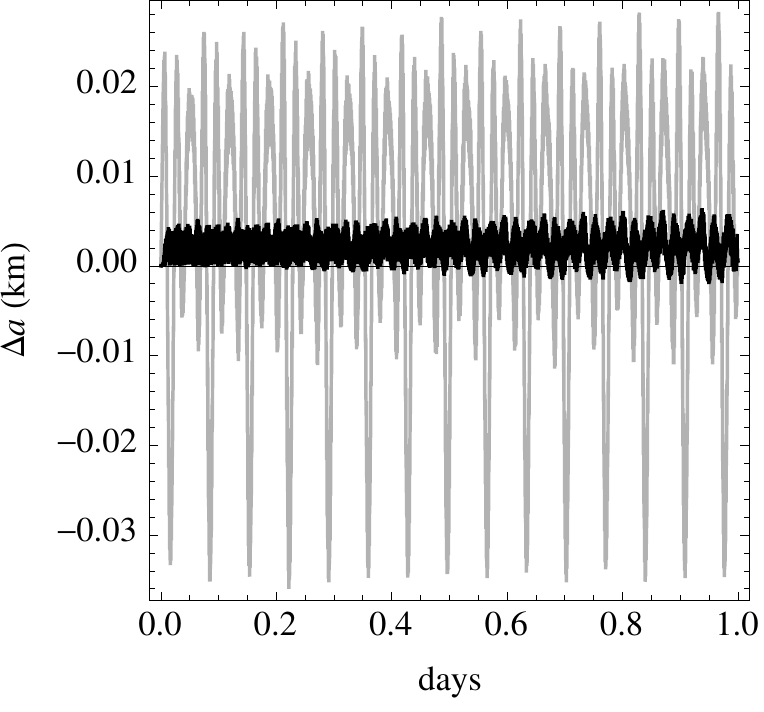}\qquad
\includegraphics[scale=0.8]{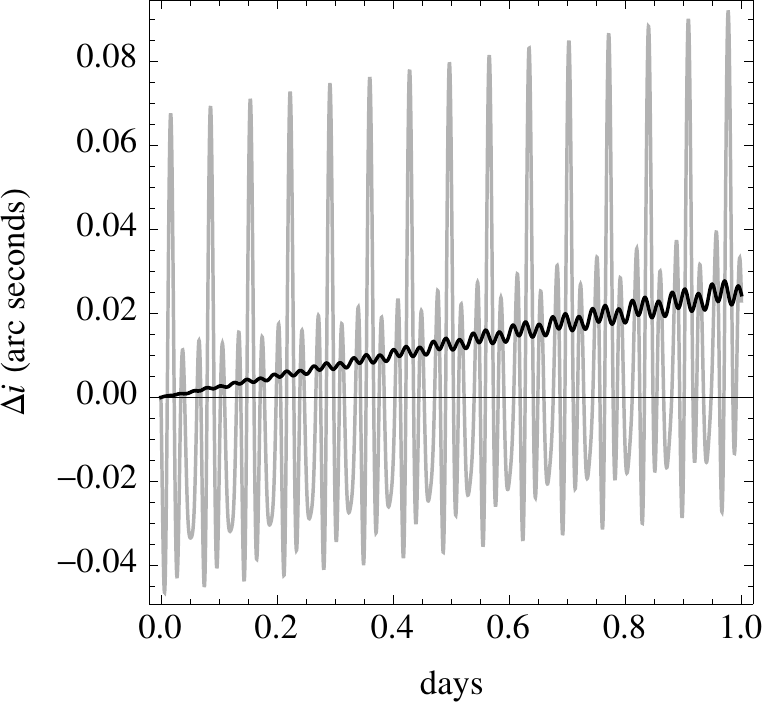}
}
\centerline{
\includegraphics[scale=0.8]{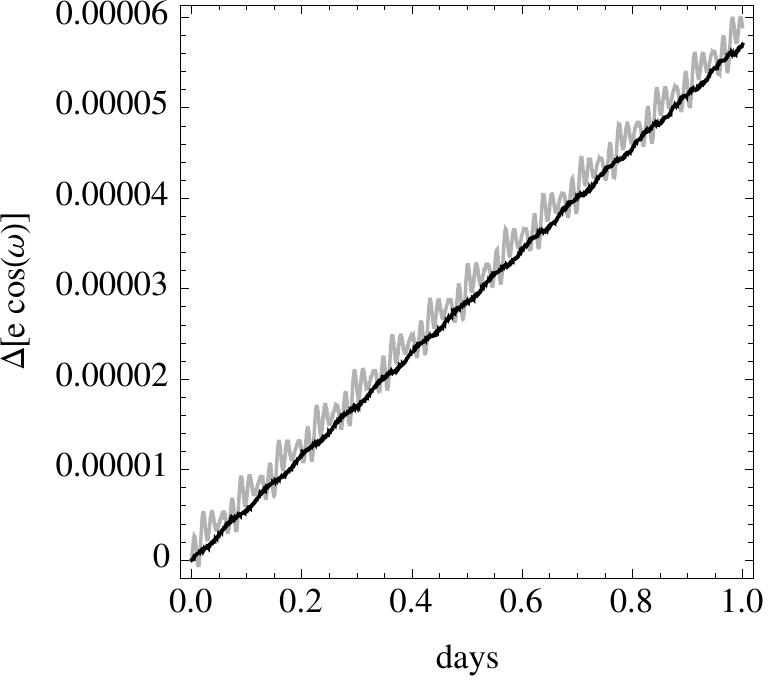}\qquad
\includegraphics[scale=0.8]{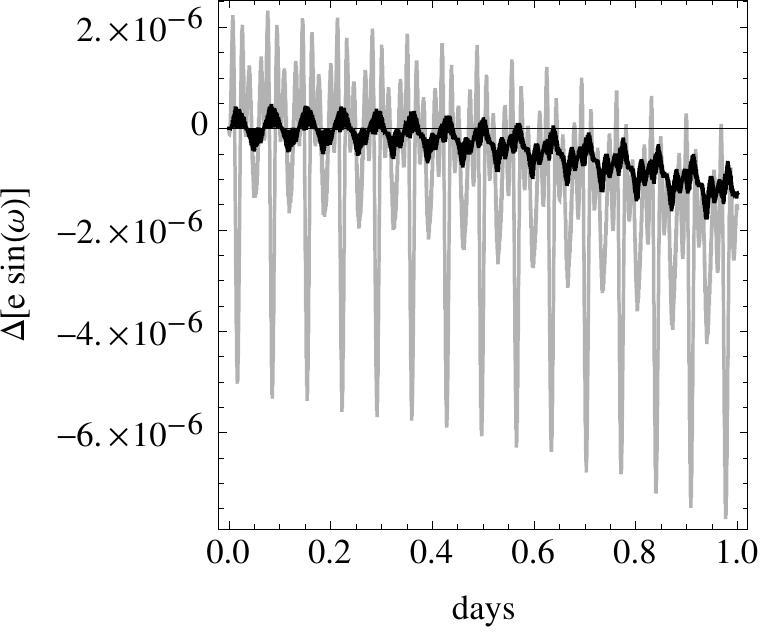} 
}
\centerline{
\includegraphics[scale=0.8]{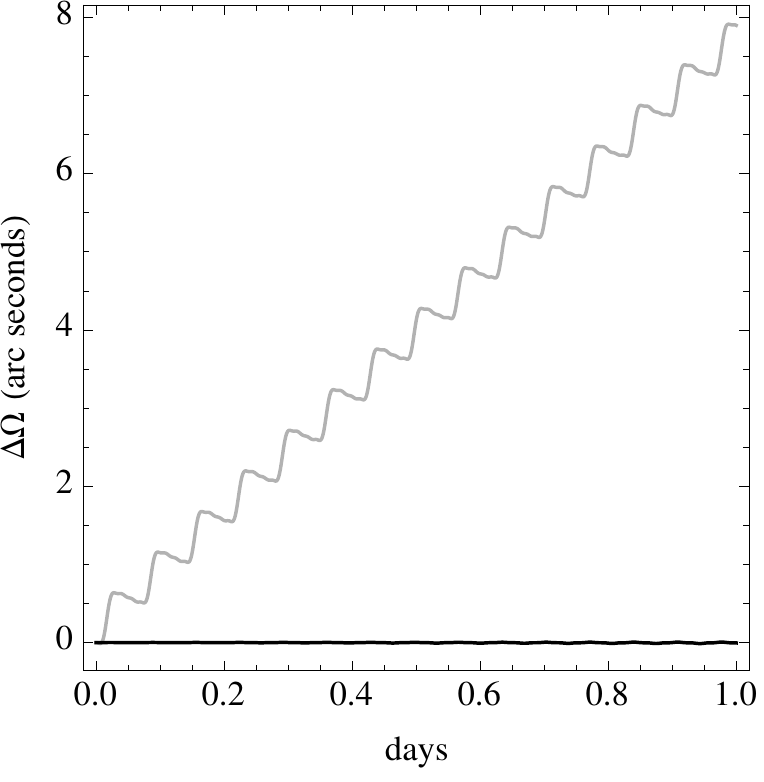}\qquad
\includegraphics[scale=0.8]{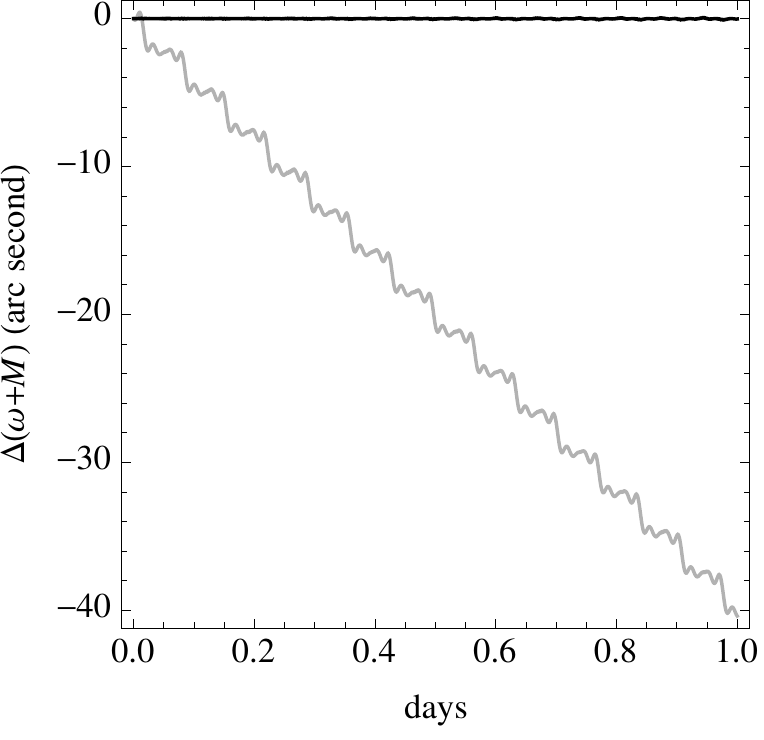}
}
\caption{Errors between the main problem and the full zonal model (gray lines) and between the {first} intermediary and the full zonal model (black lines) for a Spot-type satellite.}
\label{f:spotZvsJ2}
\end{figure}

First of all, we illustrate the effects of neglecting second order terms of the earth's gravitational potential for a Spot-type satellite. We use the initial conditions corresponding to the orbital elements given in the first row of Table \ref{t:iicc} and propagate them for 1 day (about 15 orbital periods). Errors between the full zonal model propagation and the $J_2$ truncation are shown in Fig.~\ref{f:spotZvsJ2} (gray lines), to which corresponding errors between the full zonal model propagation and the intermediary propagation (black lines) have been superimposed. Note that because of the low eccentricity of the Spot orbit, in order to avoid additional errors introduced by the inaccurate determination of the argument of perigee, instead of providing errors for the mean anomaly, the argument of the perigee, and the eccentricity, we provide errors for the mean argument of latitude
\begin{equation} \label{Mplusw}
F=M+\omega,
\end{equation}
and the semi-equinoctial elements
\begin{equation} \label{CandS}
C=e\cos\omega, \qquad S=e\sin\omega.
\end{equation}
\par

\begin{table}[tb]
\begin{center}
\begin{tabular}{@{}clllrrr@{}}
& \multicolumn{1}{c}{$a$ (km)} & \multicolumn{1}{c}{$e$} & \multicolumn{1}{c}{$I$} & \multicolumn{1}{c}{$\Omega$} & \multicolumn{1}{c}{$\omega$} & \multicolumn{1}{c}{$M$} \\
\hline
SPOT4       & $7081.1390$ & $0.0158$ & $98.0$ & $164.02\phantom{0}$ & $0.0\phantom{00}$ & $0.0\phantom{00}$ \\
typical LEO & $6831.5723$ & $0.00136$ & $51.6$ & $224.8\phantom{00}$ & $280.1\phantom{00}$ & $66.5\phantom{00}$ \\
EYE-SAT     & $7078.0$    & $0.00001$ & $98.18$ & $0.0\phantom{00}$ & $0.0\phantom{00}$ & $0.0\phantom{00}$ \\
PROBA2      & $7106.1370$ & $0.00004$ & $98.3$ & $91.364$ & $-1.423$ & $180.0\phantom{00}$ \\
JASON1      & $7254.0729$ & $0.06216$ & $66.974$ & $-74.818$ & $-241.050$ & $179.726$ \\
CRYOSAT     & $7100.4651$ & $0.00252$ & $92.029$ & $-37.185 $ & $107.492$ & $51.202$ \\
ATV         & $6586.1775$ & $0.0328$ & $51.6$ & $153.480$ & $-21.395$ & $215.240$ \\ 
Labs-Dove   & $6851.946$ & $0.0012$ & $97.326$ & $0.0\phantom{00}$ & $90.0\phantom{00}$ & $0.0\phantom{00}$ \\ 
\hline
\end{tabular}
\caption{Initial conditions used in the tests. Angular variables are in degrees \label{t:iicc}}
\end{center}
\end{table}

As shown in the top-left plot of Fig.~\ref{f:spotZvsJ2}, errors in the semi-major axis are of (short) periodic nature, and one order of magnitude smaller in the case of the intermediary than in the main problem numerical integration. The latter misses $J_4$-related secular terms, a fact that has an obvious impact on the frequencies with which the angular variables evolve, and that is manifested by the linear trend in the errors of the right ascension of the ascending node and the mean argument of latitude (bottom-left and -right plots of  Fig.~\ref{f:spotZvsJ2}, respectively) which dominates over other periodic effects missed by the main problem model. This linear trend of the errors is not observable in the intermediary propagation for this short time span, as expected from the inclusion in the intermediary solution of secular terms up to $\mathcal{O}(J_2^2)$. Finally, the missed perigee dynamics in both the intermediary and the $J_2$ truncation is clearly manifested by the long-period errors of the eccentricity and inclination ---the orbital elements related to the total angular momentum, which is the conjugate momentum to the argument of the perigee. Indeed, in the short time span encompassed by the propagation, these long-period effects result in an analogous trend in the errors of both the main problem numerical integration and the intermediary analytical solution, the former being also affected by evident short-period effects that are a consequence of the missed harmonics in the main problem model (top-right plot and center-left and -right plots in Fig.~\ref{f:spotZvsJ2}). 
\par

\begin{figure}[htbp]
\centerline{
\includegraphics[scale=0.79]{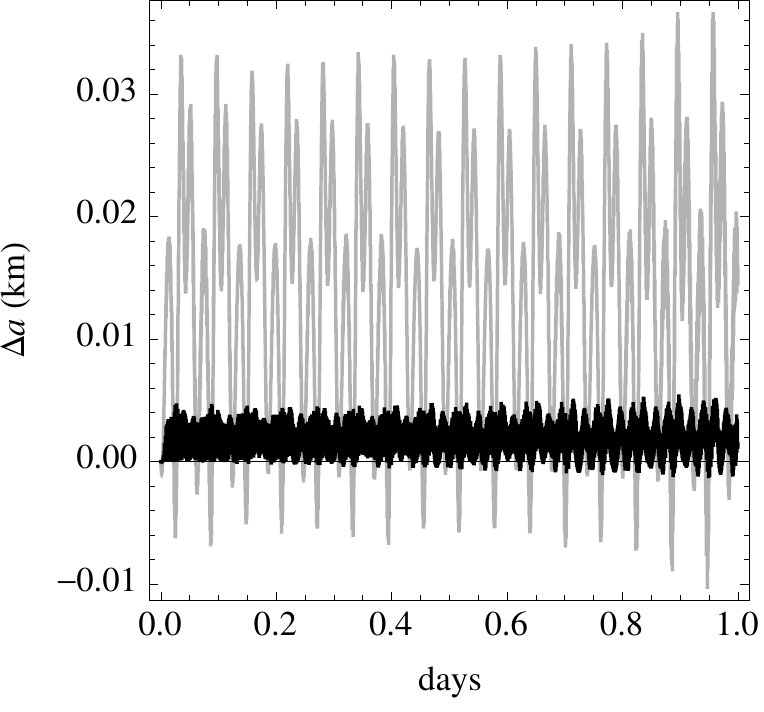}\qquad
\includegraphics[scale=0.75]{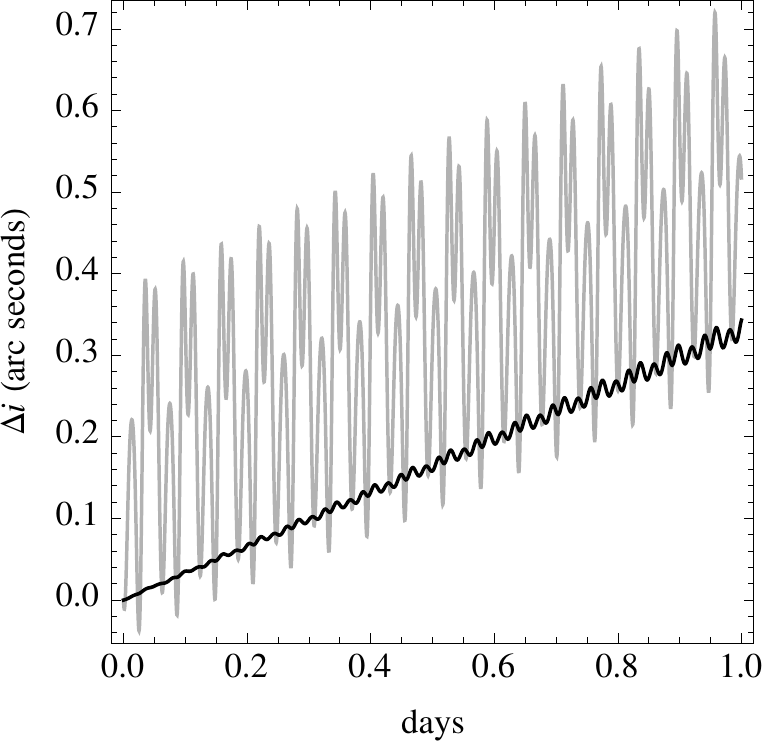}
}
\centerline{
\includegraphics[scale=0.8]{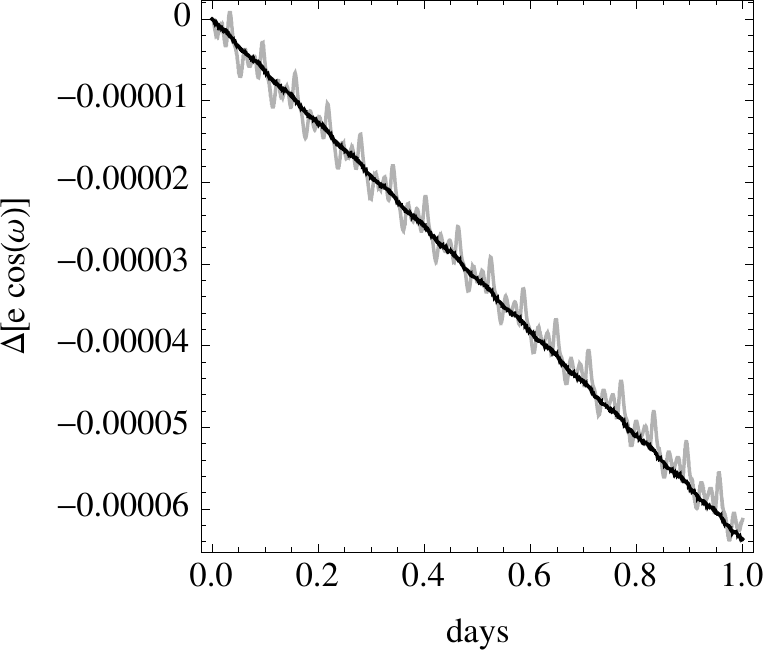}\qquad
\includegraphics[scale=0.8]{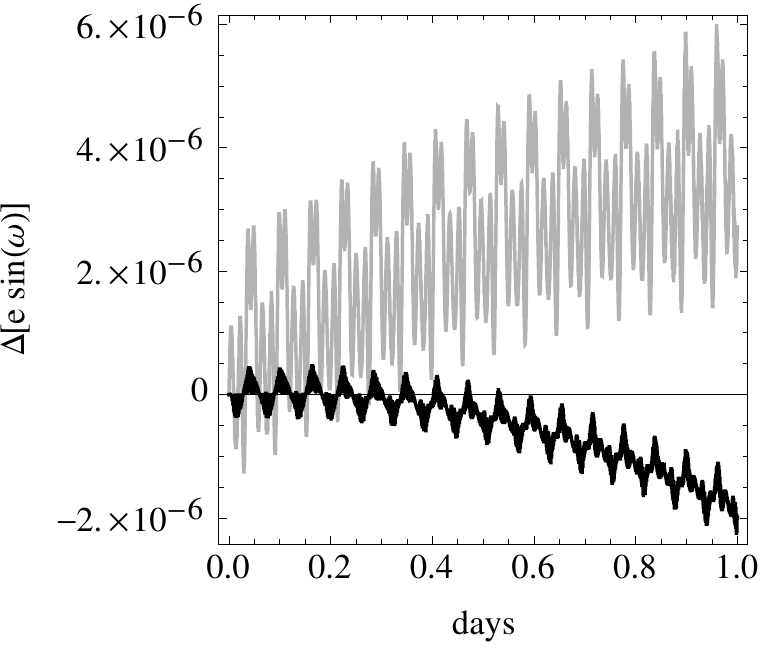} 
}
\centerline{
\includegraphics[scale=0.75]{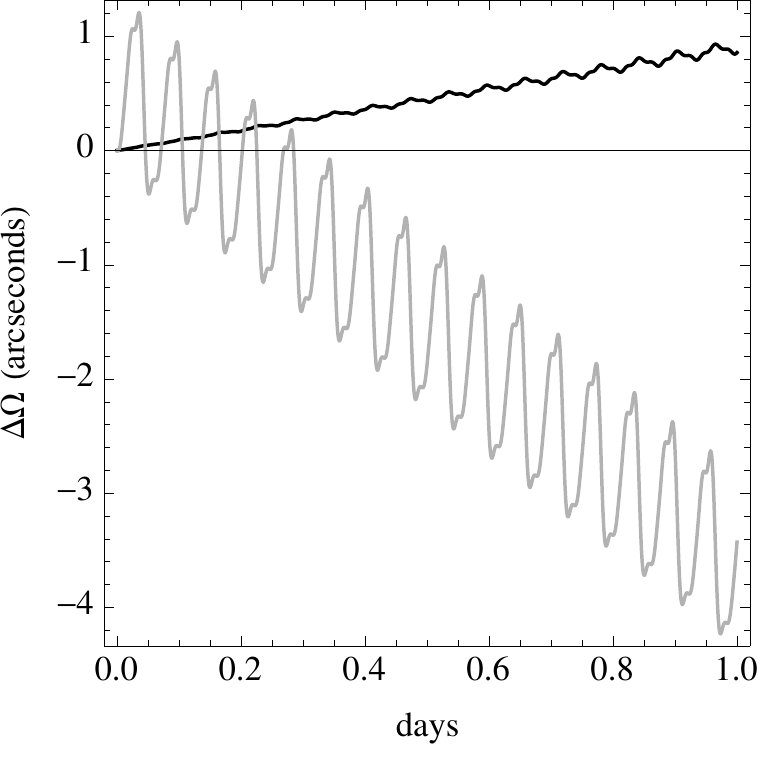}\qquad
\includegraphics[scale=0.75]{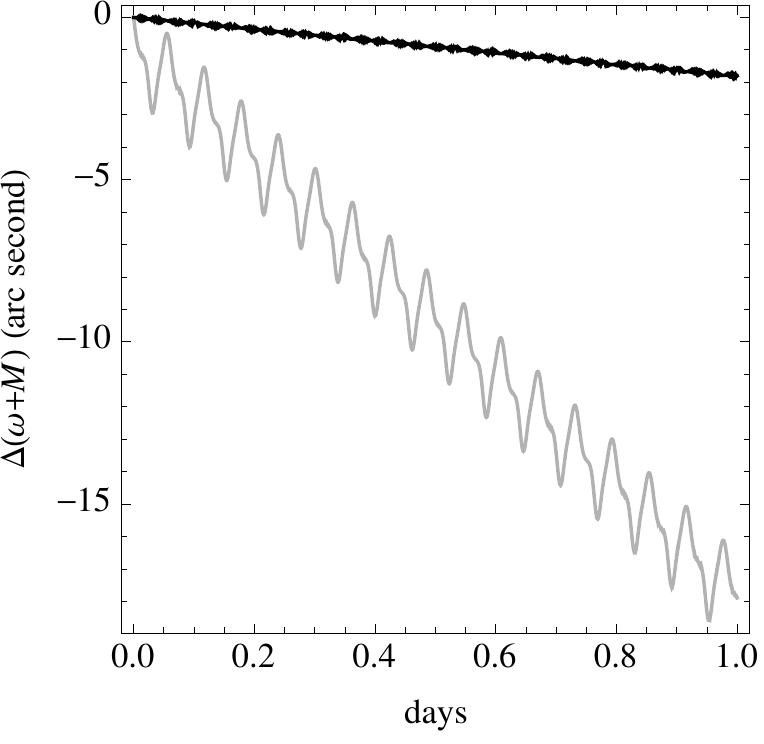}
}
\caption{Errors between the main problem and the full zonal model (gray lines) and between the {first} intermediary and the full zonal model (black lines) for one day propagation of the ATV.}
\label{f:ATV1dayJ2}
\end{figure}

Similar tests have been performed for a variety of orbits, always finding analogous results. In particular, different simulations have been carried out for the parameters of a typical LEO, but also for the nominal parameters of two Cubesats: Planet Lab-Dove and EYE-SAT, as well as the following missions: PROBA2, JASON1, CRYOSAT, and ATV. The initial conditions used are summarized in Table \ref{t:iicc}. The worst results of the intermediary propagation are obtained for JASON1 and ATV. This was expected from their higher values of the eccentricity when compared to other cases, which make errors related to the neglected long-period terms to be more apparent. Still, results provided by the intermediary are quite accurate and clearly defeat those provided by the numerical integration of the main problem, as illustrated in Fig.~\ref{f:ATV1dayJ2} for the case of ATV.

\section{Accelerating evaluation} \label{s:fe}

Keeping second order short-period effects in the computation of osculating elements may exceed the typical requirements for these kinds of simplified propagations, {as, for instance when using SGP4}. Hence, a way of speeding up evaluation of the intermediary propagation is to neglect the second order corrections of the \emph{direct} transformation (from intermediary, prime variables to osculating ones). On the contrary, in order to avoid introducing additional uncertainties in an ephemeris propagation, the initial transformation from osculating variables to intermediary variables should be as accurate as possible \citep{Lara2015ASR}. Therefore, the whole second order \emph{inverse} corrections should be taken into account. However, it must be noted that not all the orbital elements contribute in the same way to this additional uncertainty due to the truncation of the short-period, inverse transformation. Indeed, because of the known Lyapunov instability of Keplerian motion, the more relevant source of errors in an ephemeris propagation is associated to an uncertainty in the initial value of the semi-major axis, for it directly modifies the value of the mean motion ---{cf.~the thorough discussion in \citet{BreakwellVagners1970}}. How this error in the initial semi-major axis is extended to the polar-nodal elements is investigated as follows.
\par

Since the (osculating) semi-major axis is closely related to the energy of the Keplerian motion, which in polar-nodal variables writes
\begin{equation} \label{Keplerian}
\mathcal{H}_\mathrm{K}=\frac{1}{2}\left(R^2+\frac{\Theta^2}{r^2}\right)-\frac{\mu}{r},
\end{equation}
it is expected that the second order corrections of the transformation equations of the elimination of the parallax for $r$, $R$, and $\Theta$ will have the higher impact in the uncertainty with which the initial conditions are computed in the prime variables.
\par

If we replace the polar-nodal variables $r$, $R$, and $\Theta$ in Eq.~(\ref{Keplerian}), by corresponding ones given by the transformation
$\xi'=\xi+\epsilon\,\Delta_2\xi+\frac{1}{2}\epsilon^2\,\delta_2\xi$, $\xi\in(r,R,\Theta)$ {---see Appendix \ref{a:parallax} for further information on the corrections $\Delta_2$ and $\delta_2$}, it can be checked that the contribution of the second order corrections to the Keplerian energy is
\begin{equation} \label{DeltaE}
\Delta{E}=\frac{1}{2}\epsilon^2
\frac{\Theta^2}{r^2}\left[
\left(\frac{rR}{\Theta}\right)^2\frac{\delta_2R}{\Theta/p}+\frac{\delta_2\Theta}{\Theta}+\left(\frac{r}{p}-1\right)\frac{\delta_2r}{p}\right],
\end{equation}
where
\[
\frac{\delta_2R}{\Theta/p}, \qquad
\frac{\delta_2\Theta}{\Theta}, \qquad
\frac{\delta_2r}{p},
\]
are nondimensional quantities of order one. Then, in view of the coefficients $({rR}/{\Theta})$ and $(-1+{r}/{p})$ in Eq.~(\ref{DeltaE}) are $\mathcal{O}(e)$, in the case of LEO in which the present research focusses, the most important source of uncertainties in the second order corrections from original to prime elements, is associated to the second order inverse correction of the modulus of the total angular momentum.
\par

Therefore, the algorithm implementing the intermediary propagation is simplified by neglecting all the second order short-period corrections of the \emph{direct} transformation, on the one hand, and neglecting those terms of the second order \emph{inverse} transformation that contribute terms of the order of $e^2$ to Eq.~(\ref{DeltaE}), on the other hand. In particular, the only second order inverse corrections that are taken into account are: the second order correction to $\Theta$, which is taken with full precision up to $\mathcal{O}(e^2)$, and the second order correction to $r$, which is simplified by retaining only terms of $\mathcal{O}(e)$. Explicit expressions of the corrections used in the simplified version of the intermediary algorithm are given in Appendix \ref{a:parallax}. Complete second order corrections, which were thoroughly used when the intermediary was under development, can be found in \citep{HautesserresLara2016}, in which terms of the order of $e^2J_{2}^2$ and higher were neglected in agreement with other LEO simplifications.
\par

\begin{figure}[htb]
\centerline{
\includegraphics[scale=0.75]{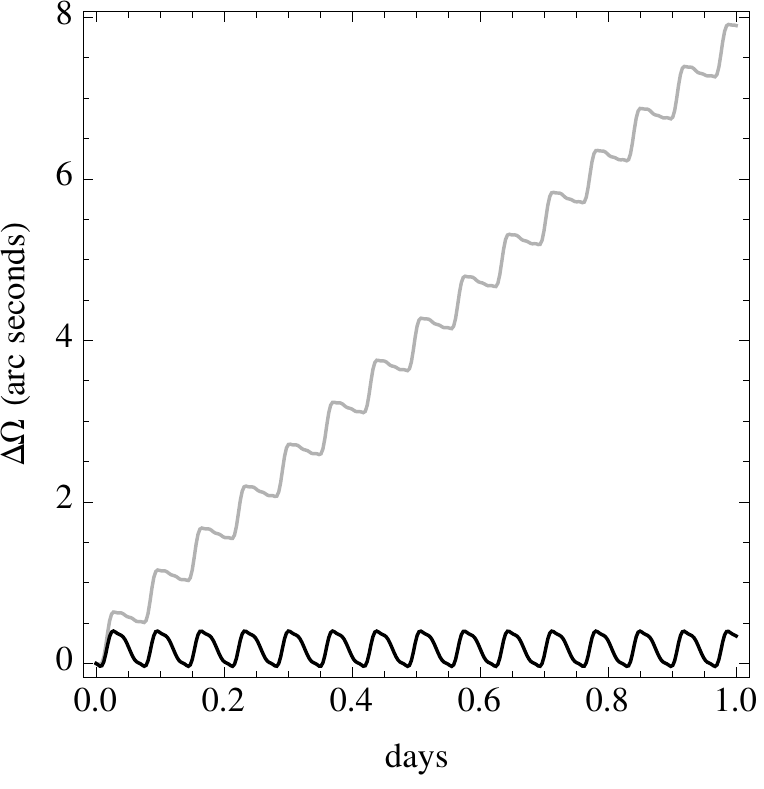}\qquad
\includegraphics[scale=0.8]{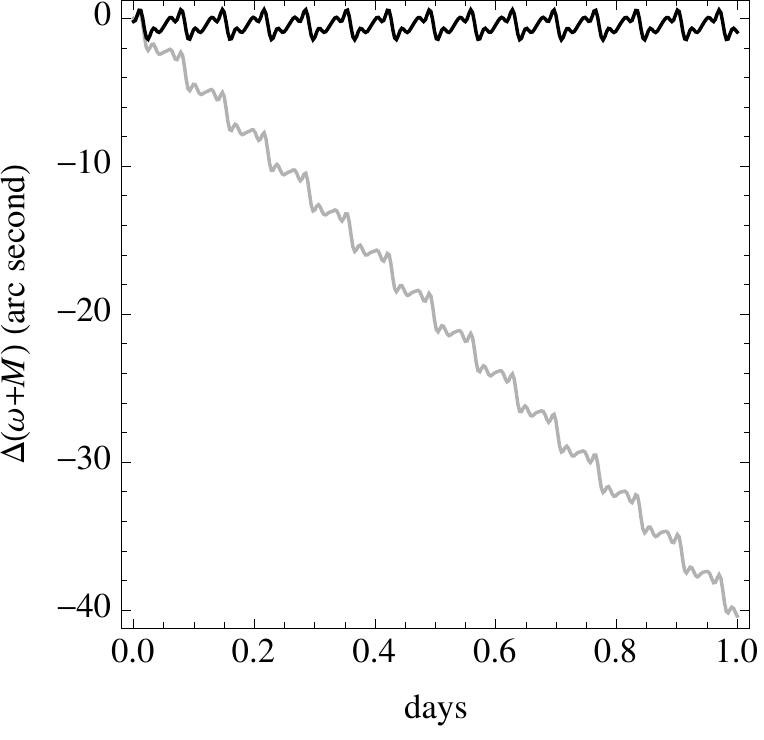} 
}
\caption{Time history of the first 24 h of the errors of $\Omega$ and $M+\omega$ for the main problem numerical integration (gray lines) and the simplified first intermediary propagation (black lines) of a Spot-type satellite.}
\label{f:spotZvsJ2fast1}
\end{figure}

{From our numerical tests, we found that} these simplifications radically reduce the effort in evaluating the intermediary to about one third of the case in which full second order corrections are used. Now, errors due to short-period terms are of comparable amplitude in both the intermediary and the $J_2$ numerical propagation, and no differences between both approaches are appreciated in the errors of the semi-major axis, inclination, and semi-equinoctial elements for one day interval. Still, the fact that the intermediary includes secular terms of the second order of $J_2$, makes that the errors of the intermediary propagation in the mean argument of latitude and right ascension of the ascending node remain free from the secular trend that affects corresponding errors of the main problem integration. This behavior is illustrated in Fig.~\ref{f:spotZvsJ2fast1} for the same test case provided in Fig.~\protect\ref{f:spotZvsJ2}.
\par

\section{Long period corrections: elimination of the perigee} \label{s:lpc}

Finding integrability of the Hamiltonian after the elimination of the parallax given in Eq.~(\ref{Kfull}), has been possible because we neglected the second order effects of the perigee dynamics. However, integrability is at the cost of missing the long-period effects associated to the odd zonal harmonics. In spite of this approximation clearly provides much better results than the usual numerical integration of the $J_2$ model, it experiences long-period errors in the time history of the eccentricity and argument of the perigee, which are apparent since the beginning of the propagation even in the case of the lower eccentricity orbits.
\par

It emerges, then, the question on how to improve the analytical solution without deteriorating evaluation performance. A plain answer is making and additional transformation: the elimination of the perigee \citep{AlfriendCoffey1984}. Thus, instead of directly ignoring the terms that are factored by the eccentricity in Eq.~(\ref{Kfull}), we only drop those terms that are multiplied by the \emph{square} of the eccentricity. Then, the only appearance of the argument of the perigee in the simplified Hamiltonian is associated to the odd zonal harmonic $J_3$. Now, this long-period term is removed by means of a canonical transformation $(\ell',g',h',L',G',H')\longrightarrow(\ell'',g'',h'',L'',G'',H'')$ based on the generating function
\begin{equation} \label{perigeeGen}
W=G\,\epsilon_3\,se\cos{g}+\mathcal{O}(e^2J_{2}),
\end{equation}
where 
\begin{equation} \label{eps3}
\epsilon_3=\frac{1}{2}\frac{J_{3}}{J_{2}}\frac{\alpha}{p}.
\end{equation}

Note that, as recommended in \citet{LaraSanJuanLopezOchoa2013c} the elimination of the perigee is carried out directly in Delaunay variables, rather than in the polar-nodal variables in which the original transformation was devised. This direct approach is straightforward and avoids dealing with the special algebra of the state functions $\tilde{C}\equiv\tilde{C}(r,\theta,R,\Theta)$ and $\tilde{S}\equiv\tilde{S}(r,\theta,R,\Theta)$, given by
\[
\tilde{C}=\frac{G}{p}C=\left(\frac{\Theta }{r}-\frac{\Theta }{p}\right)\cos\theta+R\sin\theta, \qquad
\tilde{S}=\frac{G}{p}S=\left(\frac{\Theta }{r}-\frac{\Theta }{p}\right)\sin\theta-R\cos\theta,
\]
in which canonical simplifications were originally developed \citep[cf.][]{Deprit1981}. However, computing the long-period transformation of the Delaunay variables has the inconvenience of introducing the eccentricity in denominators. Therefore, the long-period corrections of the elimination of the peri\-gee are better formulated by replacing $\ell$, $g$ and $G$ by the non-canonical variables $F$, $C$, and $S$ defined in Eq.~(\ref{Mplusw}) and (\ref{CandS}), cf.~\citep{DepritRom1970}.
\par

Thus, we obtain the long-period, first order corrections 
\begin{eqnarray} \label{dF}
\Delta{F} &=& \epsilon_3\left[\frac{1}{s}-s\left(1+\eta+\frac{1}{1+\eta}\right)\right]C, \\
\Delta{S} &=& \epsilon_3\left[C^2 \left(\frac{1}{s}-s\right)-s \left(1-S^2\right)\right], \\
\Delta{C} &=&-\epsilon_3\left(\frac{1}{s}-2s\right)CS, \\
\Delta{h} &=&-\epsilon_3\frac{c}{s}C, \\
\Delta{H} &=& 0, \\ \label{dL}
\Delta{L} &=& 0,
\end{eqnarray}
{ which in the direct transformation (resp.~inverse) must be evaluated in second prime variables (resp.~prime) and taken with the plus sign (resp.~minus). Equations (\ref{dF})--(\ref{dL})} are valid except for the lower inclinations. This is not a major issue for earth orbits in LEO, but, if required, the long period corrections can be computed in a set of non-singular variables, as, for instance, the set $(\psi,\xi,\chi,r,R,\Theta)$ of non-canonical variables \citep{Lara2015MPE}, where
\begin{equation} \label{nsnc}
\psi=\theta+\nu, \qquad
\xi=\sin{I}\sin\theta, \qquad 
\chi=\sin{I}\cos\theta.
\end{equation}
Corresponding long-period corrections are given by the sequence
\begin{eqnarray} \label{dyns}
\Delta\psi &=& \epsilon_3\left(2\chi+\frac{\kappa\chi-c\xi\sigma}{1+c}\right),
\\  \label{dxins}
\Delta\xi &=& 
\epsilon_3\left[2 \chi ^2+\kappa  \left(1-\xi ^2\right) \right],
\\ \label{dchins}
\Delta\chi &=&
-\epsilon_3\left[c^2 \sigma +(2+\kappa) \xi  \chi \right],
\\ \label{drns}
\Delta{r} &=& \epsilon_3\,\xi\,p,
\\ \label{dRRns}
\Delta{R} &=&
\epsilon_3\,(1+\kappa)\,\chi\frac{\Theta}{r},
\\ \label{dZZns}
\Delta\Theta &=&\epsilon_3(\kappa\xi-\sigma\chi)\,\Theta,
\end{eqnarray}
where $c=\sqrt{1-\xi^2-\chi^2}$, $p=\Theta^2/\mu$, and
\begin{equation} \label{kasi}
\kappa=\frac{p}{r}-1,\qquad \sigma=\frac{pR}{\Theta}.
\end{equation}
{When Eqs.~(\ref{dyns})--(\ref{dZZns}) are used for computing direct (resp.~inverse) corrections they must be evaluated in the second prime (resp.~prime) variables and added (resp.~subtracted) to the corresponding variables.}
\par

Note that $\sigma$ and $\kappa$ are $\mathcal{O}(e)$ and, therefore, Eqs.~(\ref{dyns})--(\ref{dZZns}) have no contributions of $\mathcal{O}(e^2)$. It must be noted, however, that while these corrections perform well for direct corrections, they do not in the case of inverse corrections, a case in which, as mentioned before, we should try to be as accurate as possible. 
\par

Indeed, in view of we are dealing with a zonal model, the semi-major axis is free from long-period corrections and it should not be affected by the elimination of the perigee transformation, a feature that is obviously preserved when using Eqs.~(\ref{dF})--(\ref{dL}). Quite on the contrary, when removing long-period terms using Eqs.~(\ref{dyns})--(\ref{dZZns}) the semi-major axis value is preserved only up to first order effects of $J_2$, but it is slightly modified with spurious long-period terms of $\mathcal{O}(J_2^2)$, a fact that causes an adverse impact in the computation of the initial conditions of the quasi-Keplerian system. Hence, the inverse transformation of the long-period elimination is rather computed using Eqs.~(\ref{dF})--(\ref{dL}), although it requires to compute first an additional solution of the Kepler equation, whereas computing the direct transformation using Eqs.~(\ref{dyns})--(\ref{dZZns}) is accurate enough and makes unnecessary to compute additional solutions of the Kepler equation.
\par

In addition, in order to speed computations, the inverse correction to the inclination angle can be directly computed as
\begin{equation} \label{deltaIlp}
\Delta{I}=-\epsilon_3cS.
\end{equation}
Moreover, calling $\Psi=\ell+g+h$, we compute
\[
\Delta\Psi=\Delta{F}+\Delta{h}=\epsilon_3\left(\frac{c}{1+c}+\eta+\frac{1}{1+\eta}\right)sC,
\]
which is also free from singularities in the case of low inclination orbits. Then, in the LEO assumption of low eccentricity, the inverse long-period corrections can be further simplified by neglecting $\mathcal{O}(e^2J_2)$, to obtain
\begin{equation} \label{deltaslpSimp}
\Delta\Psi=\epsilon_3\,\frac{3+5c}{2(1+c)}sC, \qquad
\Delta{S}=\epsilon_3\,s, \qquad
\Delta{C}=0,
\end{equation}
{which are evaluated in prime variables}, thus negligibly increasing the computational load of the intermediary.
\par

It is not a surprise that the long-period corrections in Eqs.~(\ref{deltaIlp}) and (\ref{deltaslpSimp}) are exactly the same as the long-period gravitational corrections used in the SGP4 propagator \citep{HootsRoehrich1980,ValladoCrawfordHujsakKelso2006}, {where they are used as direct corrections and hence have the opposite sign}. This is just a consequence of our (and theirs!) approximations for the case of low-eccentricity orbits, and the fact that the generating function of Alfriend and Coffey's elimination of the perigee applied to the $J_3$ zonal harmonic, displayed in Eq.~(\ref{perigeeGen}), matches the corresponding term of the generating function of Brouwer's long period averaging based on von Zeipel's algorithm \citep{Brouwer1959}. This feature, which {becomes clearly apparent} when the elimination of the perigee is carried out in Delaunay variables, does not seem to have been previously noticed ---probably because the traditional way in which the elimination of the perigee is computed provides the generating function in the form of a mixture of polar nodal variables and parallactic functions that complicates things unnecessarily \citep[see][for further details on this topic]{LaraSanJuanLopezOchoa2013c}. 
\par

\begin{figure}[htbp]
\centerline{
\includegraphics[scale=0.8]{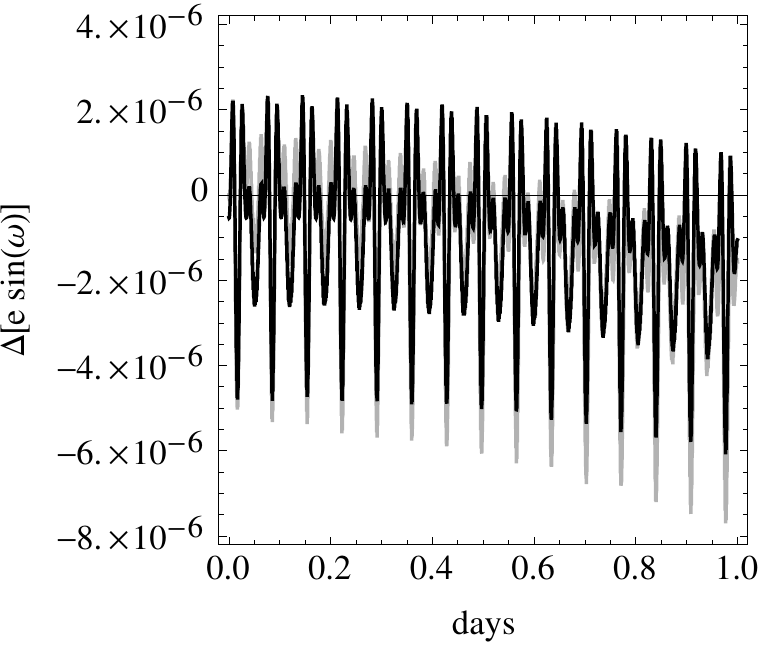}\qquad
\includegraphics[scale=0.8]{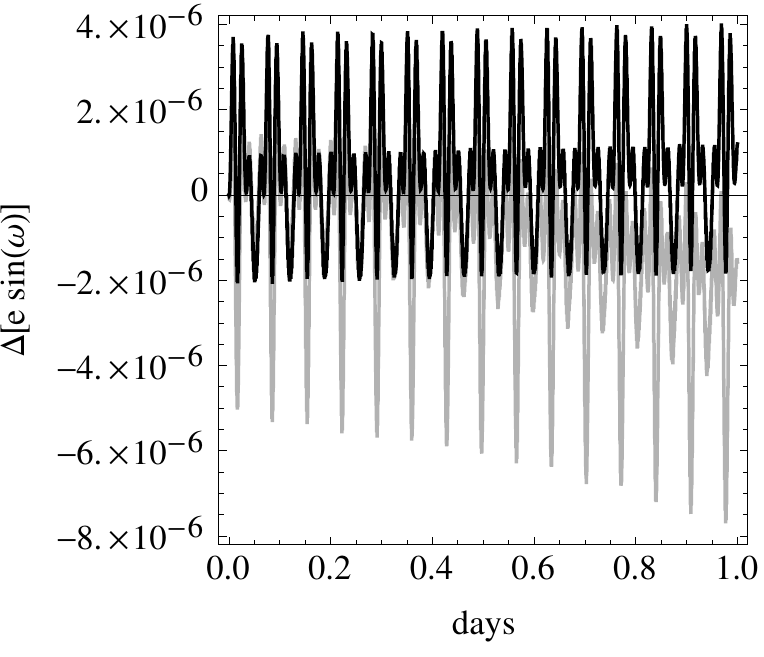}
}
\centerline{
\includegraphics[scale=0.8]{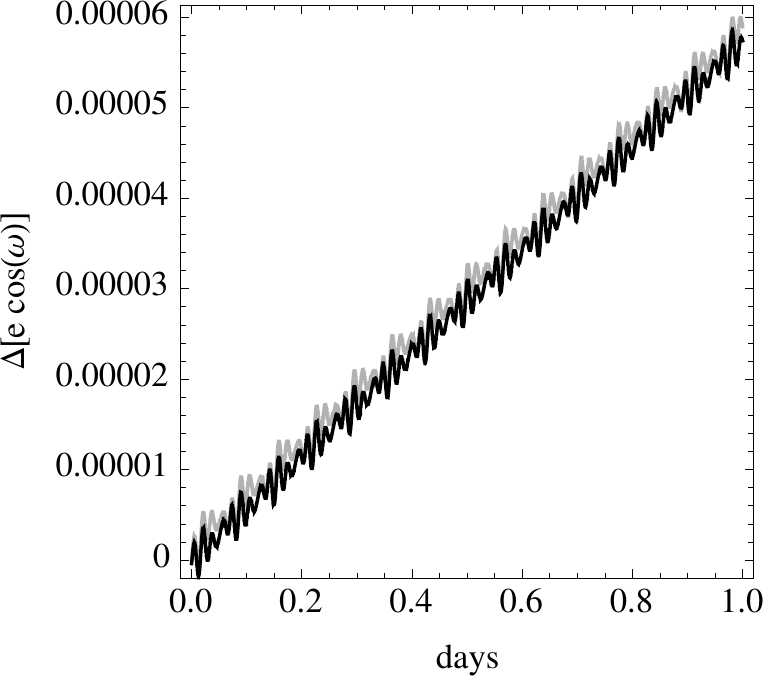}\qquad
\includegraphics[scale=0.8]{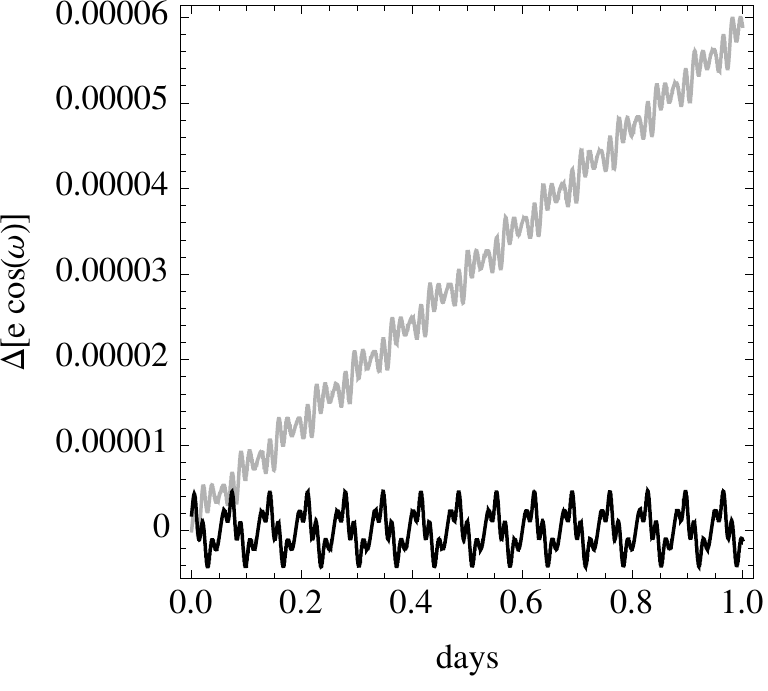} 
}
\centerline{
\includegraphics[scale=0.8]{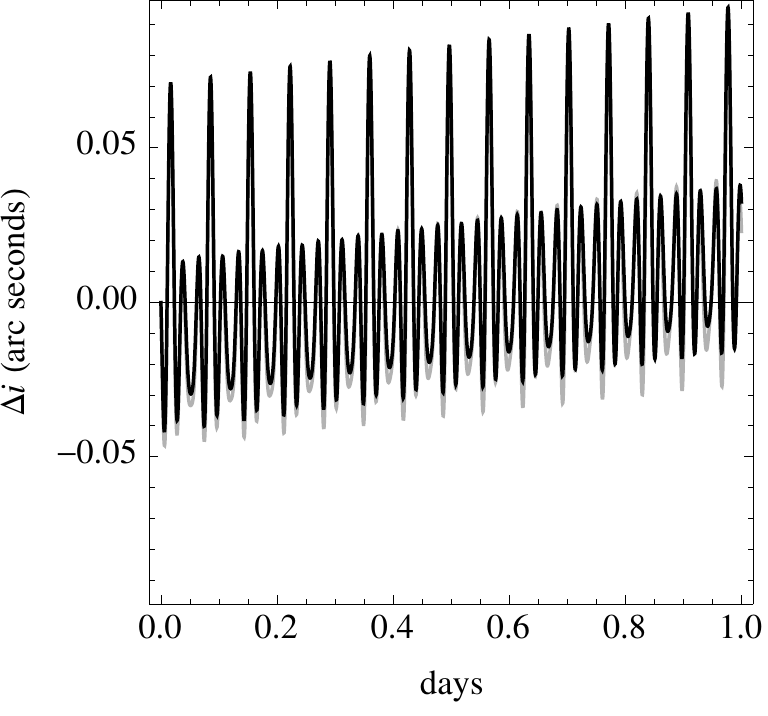}\qquad
\includegraphics[scale=0.8]{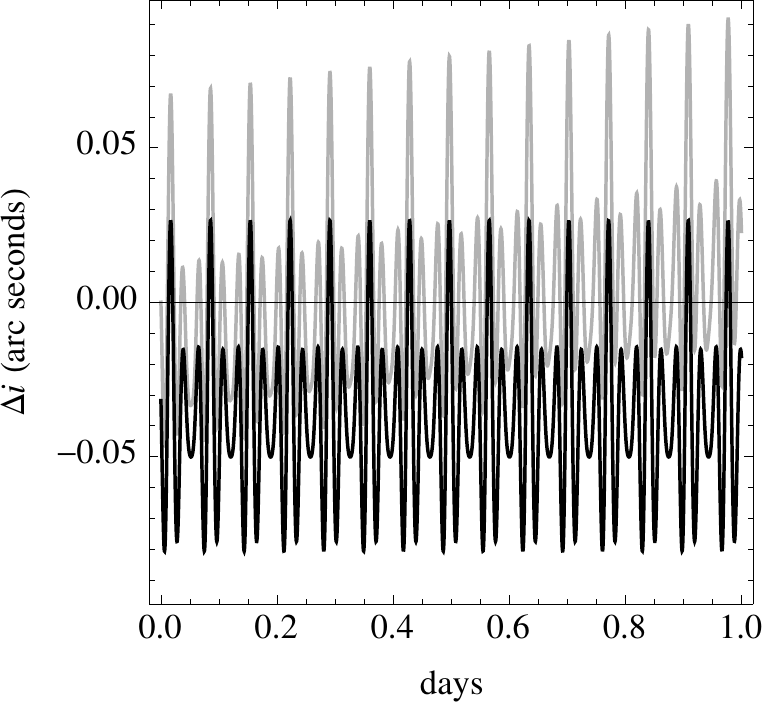}
}
\caption{One day propagation of the Spot satellite: errors of the accelerated intermediary (black lines) without (left column) and with (right column) long-period corrections, superimposed to the errors of the main problem numerical propagation (gray lines).}
\label{f:spotZvsJ2lp}
\end{figure}

The improvements achieved when the intermediary solution is upgraded with long-period corrections are illustrated in Fig.~\ref{f:spotZvsJ2lp} for the Spot-type satellite example. Plots in the left column of Fig.~\ref{f:spotZvsJ2lp} correspond to the accelerated intermediary without long-period corrections and those in the right column do incorporate the long-period corrections to the accelerated intermediary. The removal of most of the long-period errors from the intermediary propagation is evident now from a simple comparison of the left and right columns of Fig.~\ref{f:spotZvsJ2lp}.
\par

\begin{figure}[htbp]
\centerline{
\includegraphics[scale=0.8]{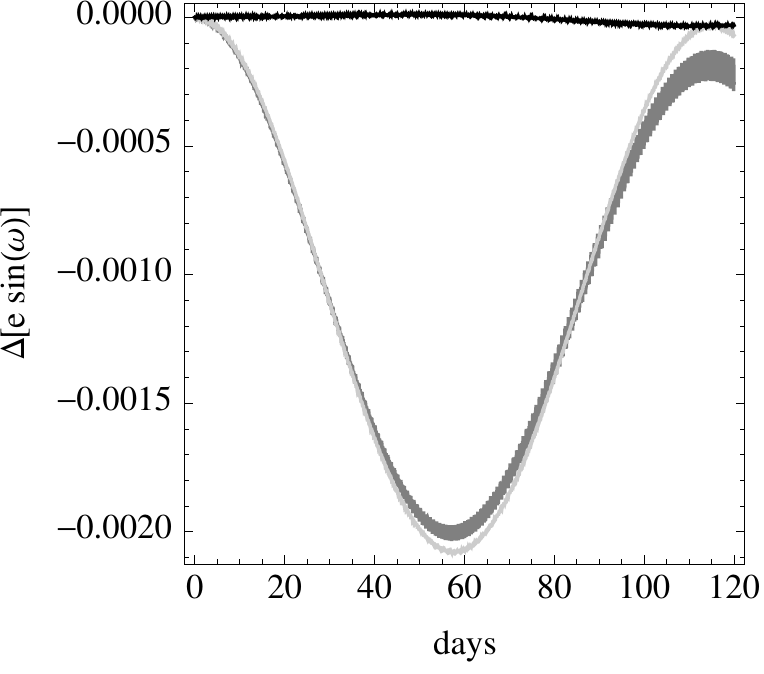}
\includegraphics[scale=0.8]{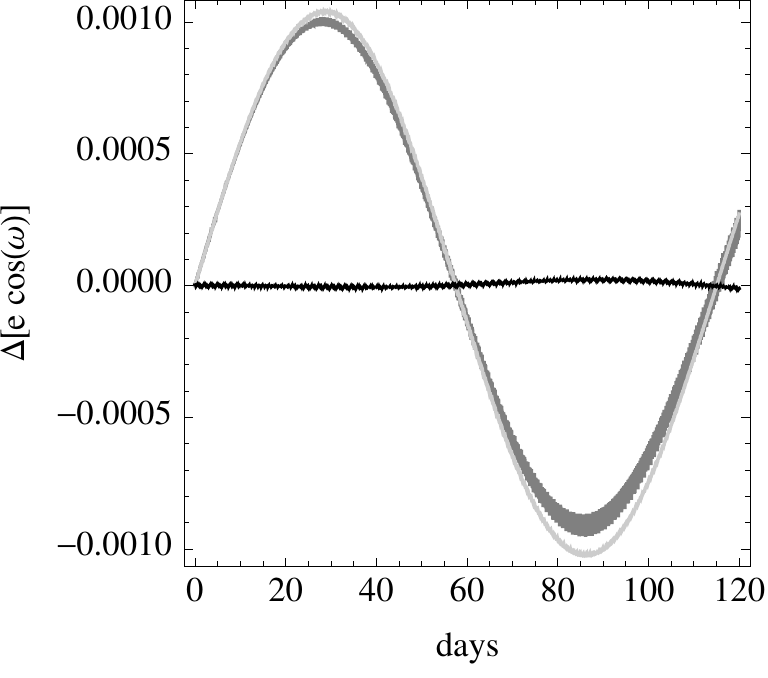} 
}
\centerline{
\includegraphics[scale=0.8]{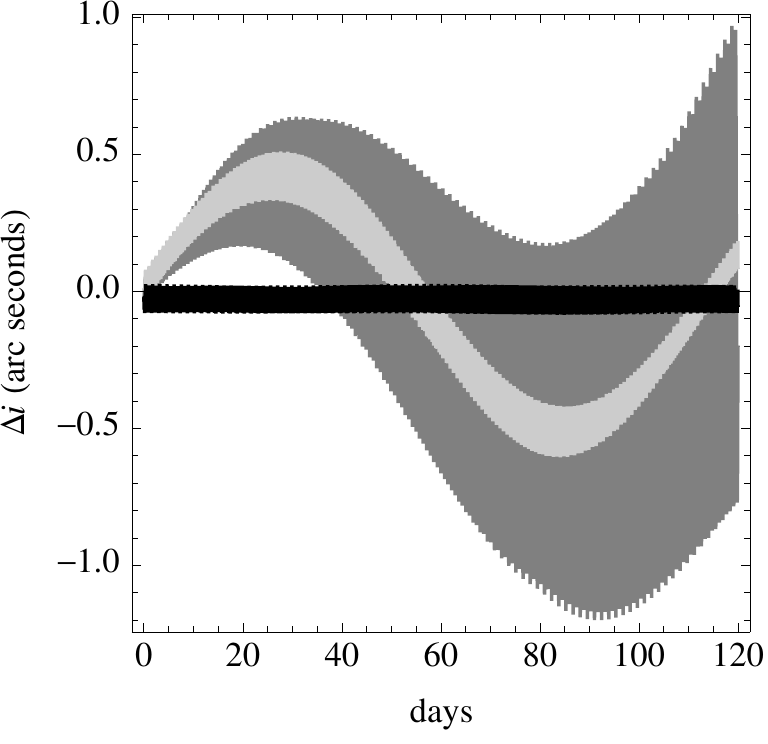}
}
\caption{Four months propagation of a Spot-type satellite. Superimposed to the errors of the main problem numerical propagation (dark gray lines) are those of the accelerated intermediary without (dark light lines) and with long-period corrections (black lines).}
\label{f:spotZvsJ2lp4}
\end{figure}

The radical improvement of the intermediary performance obtained when long-period corrections are taken into account is better illustrated when the propagation is extended to four months, a time interval that encompasses one full period of the perigee motion of the Spot satellite. Corresponding results are presented in Fig.~\ref{f:spotZvsJ2lp4}, where the effects of neglecting the long-period terms associated to the $J_3$ harmonic coefficient are clearly apparent. In particular, a long-period modulation of the errors of the intermediary propagation with the same period as the argument of the perigee of the Spot satellite (of about 16 weeks and a half) is now clearly noted in the inclination and the eccentricity. Because of that, the errors of the elements $C$ and $S$ are very similar in the $J_2$ numerical integration and in the intermediary evaluation when long-period effects are neglected, although for the later the amplitude of the short-period errors is almost negligible because of the first order short-period corrections used by the intermediary. Quite on the contrary, the long-period errors of the intermediary propagation reduce to a minimum when the elimination of the perigee is included in the algorithm.
\par

\section{Summary and runtime comparisons} \label{s:recapitulation}

A cascade of canonical transformations $\mathcal{T}_i$ $(i=1,2,3)$, allowed us to reduce, after truncation, the 2-DOF zonal Hamiltonian $\mathcal{H}\equiv\mathcal{H}(\ell,g,-,L,G,H)$ to a pure Keplerian system (in new variables):
\[
\mathcal{T}_3\circ\mathcal{T}_2\circ\mathcal{T}_1:\mathcal{H}(\ell,g,-,L,G,H)\longrightarrow-\frac{\mu^2}{2\tilde{L}^2}.
\] 
The sequence of symbolic manipulations that lead to the integrable solution is as follows:
\begin{itemize}
\item The elimination of the parallax $\mathcal{T}_1$ is applied to remove non-essential short-period terms up to the second order of $J_2$. The simplified Hamiltonian remains of 2-DOF, and is further simplified in the case of low eccentricity orbits:
\[
\mathcal{T}_1:\mathcal{H}(\ell,g,-,L,G,H)\longrightarrow\mathcal{H}_1(\ell',g',-,L',G',H)+\mathcal{O}(e^2J_2^2).
\]
\item Next, the elimination of the perigee $\mathcal{T}_2$ is applied to $\mathcal{H}_1$
\[
\mathcal{T}_2:\mathcal{H}_1(\ell',g',-,L',G',H)\longrightarrow\mathcal{H}_2(\ell'',-,-,L'',G'',H)+\mathcal{O}(e^2J_2^2).
\]
After truncation, the new Hamiltonian is of 1-DOF, and, therefore, integrable.
\item A torsion transformation $\mathcal{T}_3$ is then applied to $\mathcal{H}_2$ after reformulating it in polar nodal variables $\mathcal{H}_2\equiv\mathcal{K}(r,-,-,R,\Theta,N)$,
\[
\mathcal{T}_3:\mathcal{K}(r,-,-,R,\Theta,N)\longrightarrow\mathcal{Q}(\tilde{r},-,-,\tilde{R},\tilde{\Theta},\tilde{N})
\]
The final Hamiltonian $\mathcal{Q}\equiv \frac{1}{2}(\tilde{R}^2+\tilde\Theta^2/\tilde{r}^2)-{\mu}/{\tilde{r}}$ is a pure Keplerian system whose solution is standard.
\end{itemize}
This sequence gives rise to what we called the ``second'' intermediary.
\par

For the case of the lower eccentricity orbits, the elimination of the parallax can be written as
\[
\mathcal{T}_1:\mathcal{H}(\ell,g,-,L,G,H)\longrightarrow\mathcal{H}_2(\ell',-,-,L',G',H)+\mathcal{O}(e\,J_2^2)
\]
which, after truncation, can be directly converted into a pure Keplerian system by means of the torsion $\mathcal{T}_3$, giving rise to what we called the ``first'' intermediary.
\par

The sequence of computations required by the analytical propagation of a given set of initial conditions with any of the two versions of the (accelerated) intermediary solutions are summarized in Algorithm \ref{alg1}. The first intermediary is faster, but misses the second order effects of the perigee dynamics ---a fact that is clearly appreciated when the propagation errors are depicted in the $(e\cos\omega,e\sin\omega)$ chart, cf.~Figs.~\ref{f:spotZvsJ2} and \ref{f:ATV1dayJ2}. However, the accumulation of these long-period errors is small in short time intervals, a case in which the first intermediary provides more than acceptable results. On the other hand, the second intermediary deals effectively with the perigee dynamics, and hence, in spite of its higher computational load, it should be preferred for longer orbit propagation intervals.
\par

\begin{algorithm}[htbp]
\caption{\label{alg1} First and second accelerated intermediaries}
\begin{algorithmic}[1] 
\State  \textbf{Inputs}: Initial epoch $t_0$, final epoch $T$, evaluation interval $\Delta{t}$; initial conditions $(r_0,\theta_0,\nu_0,R_0,\Theta_0,N_0)$; physical parameters $\mu$, $\alpha$, $J_2$, $J_3$, $J_4$.
\State Compute prime variables from the inverse transformation Eq.~(\ref{inverse})
\Statex \quad use (with minus sign) the first order corrections in Eqs.~(\ref{Dr})--(\ref{DN})
\Statex \quad use the second order corrections in Eqs.~(\ref{d2r}) and (\ref{d2Z})
\If{second intermediary}
\State transform polar variables into orbital elements (involves solving Kepler equation)
\State compute inverse long-period corrections from Eqs.~(\ref{deltaIlp}) and (\ref{deltaslpSimp})
\State transform the corrected variables into (double prime) polar-nodal variables
\EndIf
\State Compute tilde variables:
\Statex \quad evaluate $\Phi$, ${\partial\Phi^2}/{\partial{c}}$, and ${\partial\Phi^2}/{\partial\epsilon}$ from Eqs.~(\ref{Phi2}), (\ref{dPhidc}), and (\ref{dPhidep})
\Statex \quad make $\tilde{r}=r$, $\tilde{R}=R$, $\tilde{N}=N$, and solve $\tilde\theta$, $\tilde\nu$, and $\tilde\Theta$ from Eq.~(\ref{transform})
\State Compute the constants $a,e,i,\Omega,\omega,M(t_0),n$ of the Keplerian problem (tilde variables)
\While{$t_i$ less than $T$ }
\State make $t_{i+1}=t_{i}+\Delta{t}$ and compute $M(t_{i+1})=M(t_0)+nt_{i+1}$
\State transform orbital elements into polar-nodal (tilde) variables
\State compute non-tilde variables from Eq.~(\ref{transform}). {Use Eq.~(\ref{ThetaNR}) for computing $\Theta$}
\If{second intermediary}
\State compute non-singular variables from Eq.~(\ref{nsnc})
\State apply long-period corrections using Eqs.~(\ref{dyns})--(\ref{dZZns})
\State transform the corrected non-singular variables into polar-nodal (prime) variables
\EndIf
\State recover (1st order) short-period terms: use Eqs.~(\ref{direct}), (\ref{Dr})--(\ref{DN}) with the plus sign
\State \textbf{Output}: time $t_{i+1}$ and state vector $(r(t_{i+1}),\theta(t_{i+1}),\nu(t_{i+1}),R(t_{i+1}),\Theta(t_{i+1}),N_0)$
\EndWhile
\end{algorithmic}
\end{algorithm}

In both cases the accuracy of the intermediary propagation is better than the usual Cowell propagation of the $J_2$ problem, which misses the contribution of the higher order harmonics. The intermediary propagation can also be advantageous in terms of computing time. However, because numerical and analytical methods are essentially different in nature, a direct comparison between both approaches only makes sense when the needs for ephemeris representation are clearly understood. Indeed, while analytical solutions can be directly evaluated at a given time, the numerical procedure requires a step by step integration.
\par

Hence, one must establish the particular scenario in which the method is applied and limit the conclusions to this particular scenario. In our case, we are interested in applications to onboard orbit propagation. Therefore, we assume the same scenario as in \cite{GurfilLara2014CeMDA}, in which the Cowell integration is carried out with a fourth-order R-K integration using a constant step-size of one second. 
\par

\begin{figure}[htbp]
\centering
\includegraphics[scale=0.75]{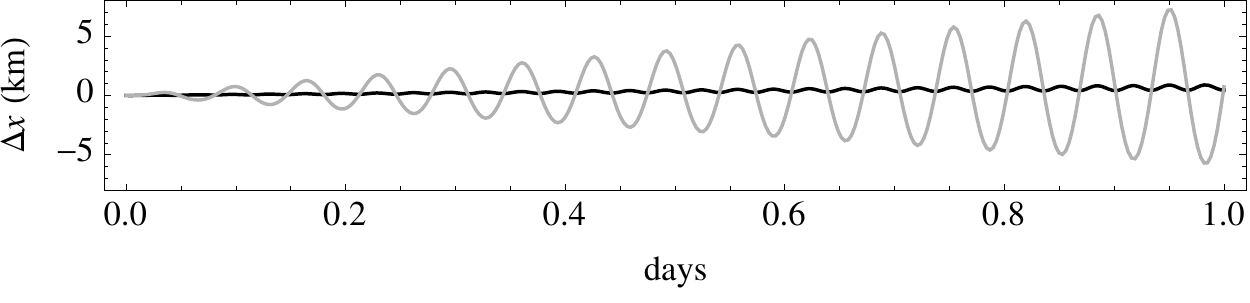}
\includegraphics[scale=0.75]{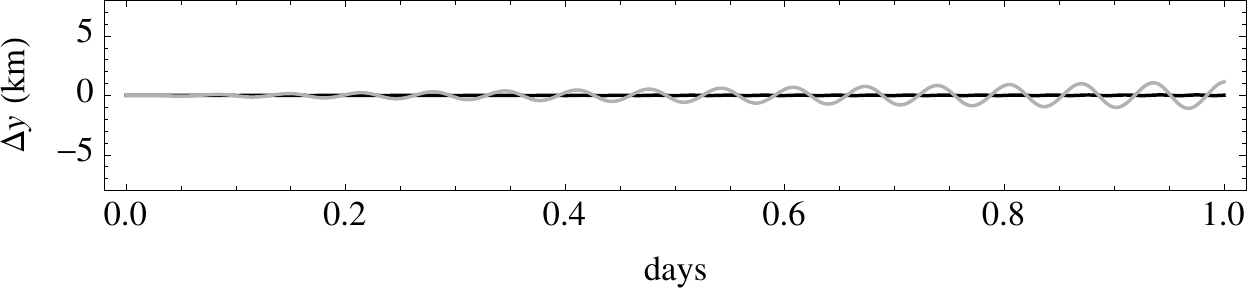}
\includegraphics[scale=0.75]{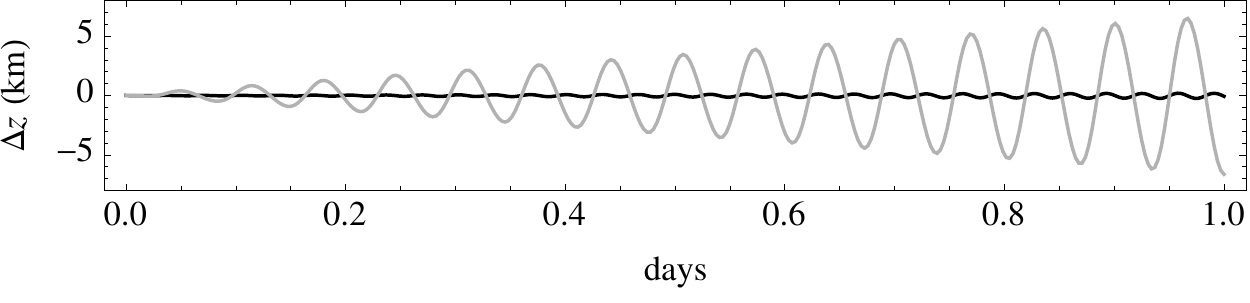}
\includegraphics[scale=0.75]{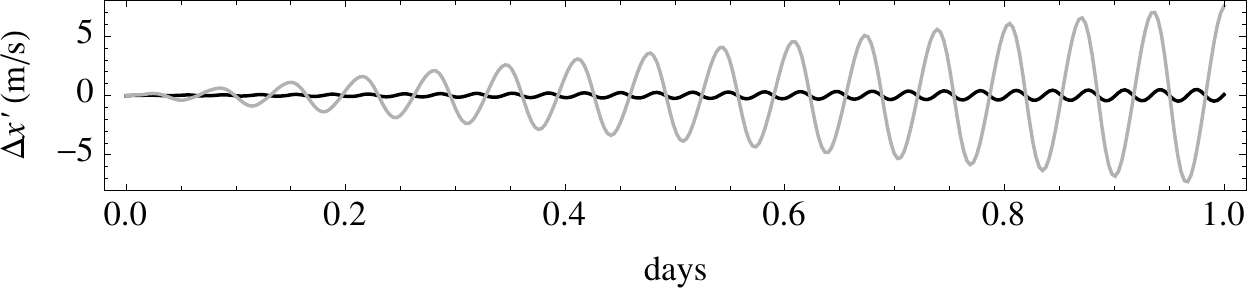}
\includegraphics[scale=0.75]{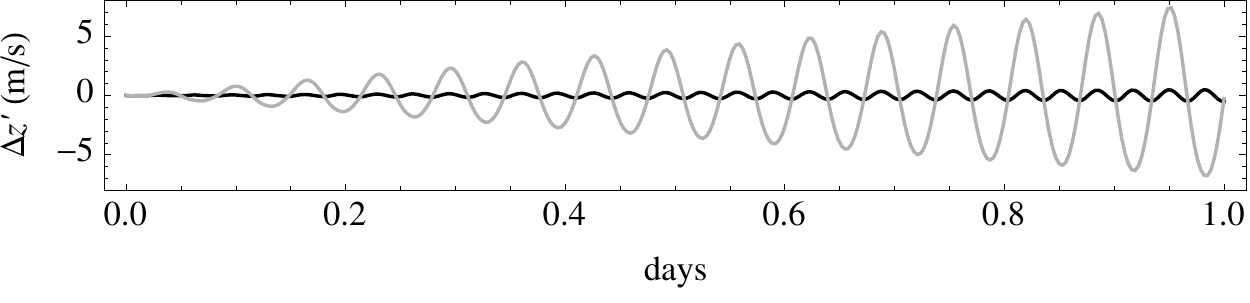}
\includegraphics[scale=0.75]{dovew}
\caption{One day propagation of a Planet Labs-Dove satellite. Superimposed to the errors of the main problem numerical integration (gray lines) are those of the accelerated first intermediary (black lines).}
\label{f:planetlab1st}
\end{figure}

To illustrate runtime comparisons, we focus on the one-day propagation of the initial conditions of a typical satellite of the Planet Labs-Dove constellation ---an almost circular sun-synchronous orbit at about 475 km of altitude--- although similar results are obtained when performing analogous tests for the other satellites in Table \ref{t:iicc}. The propagation errors in Cartesian coordinates for this case are presented in Fig.~\ref{f:planetlab1st}, where it is shown that the accuracy of the accelerated first intermediary is roughly one order of magnitude better than the $J_2$ problem R-K propagation.
\par

We checked that both, the R-K and the analytical solution take the same time when the accelerated first intermediary is evaluated each second and a half, which means a flight of the ``dove'' of about 12 km. Hence, the first intermediary will be advantageous both in terms of accuracy and computing time when ephemeris are acceptable for larger intervals than 12 km. For instance, if Lab-Dove's ephemeris are recorded each four minutes ---or more precisely, when the one day propagation interval is divided into 333 steps, which allowed to display smooth graphics in Fig.~\ref{f:planetlab1st}--- the first intermediary is about 130 times faster than the R-K integration of the $J_2$ problem.
\par

Runtime performance slightly deteriorates in the case of the second intermediary. Nevertheless, the four minutes evaluation case shows that the second intermediary still performs 100 times faster than R-K, whereas the runtime performance of both methods balance when the second intermediary is evaluated about each 2.5 seconds ---in which the dove flies $\sim20$ km.

These ratios may still improve by making additional efforts in optimizing the evaluation of the analytical solutions. In particular, Algorithm \ref{alg1} clearly discloses the \emph{profile} \citep{Knuth1970} of the intermediary procedure, which suggests that further work should concentrate in speeding up the evaluation of the different direct transformations of the intermediary solution.

\section{Conclusions}

Inclusion of higher order harmonics of the gravitational potential clearly improves the propagation model for LEO orbits, but notably penalizes the usual Cowell integration in terms of computing time. This increase of the computational load, which is directly reflected in power consumption, can be radically alleviated in the case of the lower eccentricity orbits, which are the majority in a LEO catalogue. Indeed, for the usual eccentricities of operational satellites in LEO, the Cowell propagation can be replaced by an analytical intermediary solution within a reasonable accuracy. The intermediaries proposed in the present research include higher order secular and periodic effects, which admit a compact form of straightforward evaluation when using polar-nodal variables, and hence may be adequate for onboard orbit propagation in such satellite missions in which reduced power consumption is a constraint.

\appendix

\section{Appendix}
\label{a:parallax}

Calling $\xi$ to any of the canonical variables and using
\begin{equation} \label{aepsilon}
\epsilon=-\frac{1}{2}\frac{\alpha^2}{p^2}J_{2},
\end{equation}
as given in Eq.~(\ref{epsilon}), under the simplifications used by the LEO intermediary, the direct transformation is written
\begin{equation} \label{direct}
\xi=\xi'+\epsilon\,\Delta_1\xi'
\end{equation}
and the inverse transformation
\begin{equation} \label{inverse}
\xi'=\xi+\epsilon\,\Delta_2\xi+\frac{1}{2}\epsilon^2\,\delta_2\xi.
\end{equation}

The corrections are conveniently expressed in polar nodal variables and are given below, where the eccentricity functions $\kappa$ and $\sigma$ defined in Eq.~(\ref{kasi}) are used for convenience. Besides, the abbreviation $\tilde{J}_{m}=J_{m}/J_{2}^2=\mathcal{O}(1)$, $m>2$, is used in the second order corrections.
\par

Note that the first order corrections were first provided by \cite{Deprit1981}, and are given here for the sake of completeness ---but in the arrangement proposed by \cite{GurfilLara2014CeMDA} whose evaluation is much more efficient. Besides, second order terms were previously given in \citep{Lara2015ASR} but limited to the $J_2$ perturbation. 

\subsection{First order corrections}

The first-order corrections are formally the same both in the direct and inverse transformations but with different signs, namely $\Delta_1=\Delta$ and $\Delta_2=-\Delta$ where:
\begin{eqnarray} \label{Dr}
\Delta{r} &=& p\left(1-\frac{3}{2}s^2-\frac{1}{2}s^2\cos2\theta\right), 
\\ \label{Dz}
\Delta\theta &=& \Big[1-6c^2+(1-2c^2)\cos2\theta\Big]\sigma-\left[\frac{1}{4}-\frac{7}{4}c^2+(1-3c^2)\,\kappa\right]\sin2\theta,
\\ \label{Dv}
\Delta\nu &=& c\left[(3+\cos2\theta)\,\sigma-\left(\frac{3}{2}+2\kappa\right) \sin2\theta\right],
\\ \label{DR}
\Delta{R} &=& \frac{\Theta}{r}\,(1+\kappa)s^2\sin2\theta,
\\ \label{DZ}
\Delta\Theta &=& -\Theta\,s^2\left[\left(\frac{3}{2}+2\kappa\right)\cos2\theta+\sigma\sin2\theta\right],
\\ \label{DN}
\Delta{N} &=& 0,
\end{eqnarray}
and the right member of each of Eqs.~(\ref{Dr})--(\ref{DN}) as well as $p$ in Eq.~(\ref{aepsilon}) must be expressed in prime variables when computing $\Delta_1\xi'\equiv\Delta\xi'$, or in original ones when computing $\Delta_2\xi\equiv-\Delta\xi$.

\subsection{Second order, simplified inverse corrections}

The necessary corrections are given in following formulas, where all the symbols are functions of the original variables. The correction $\delta_2{R}$ is not used by the accelerated intermediary and is provided just for convenience of those interested in checking Eq.~(\ref{DeltaE}).
\begin{eqnarray} \label{d2r}
\delta_2{r} &=& p\,\Big\{-3+10c^2+c^4 -(4-32c^2)s^2\cos2\theta-s^4\cos4\theta \\ \nonumber
&& -\frac{3}{2}\frac{p}{\alpha}\tilde{J}_{3}\Big[(1-5c^2)s\sin\theta+\frac{5}{6}s^3\sin3\theta\Big] \\ \nonumber
&& -\tilde{J}_{4}\Big[ \frac{9}{8}(3-30c^2+35c^4)+\frac{5}{2}(1-7c^2)s^2\cos2\theta-\frac{7}{8}s^4\cos4\theta \Big]
+\mathcal{O}(e)\Big\}
\\[1ex]
\delta_2{R} &=& \frac{\Theta}{p}\Big\{s^2(2-22c^2)\sin2\theta +s^4\sin4\theta
+\frac{3}{2}\frac{p}{\alpha}\tilde{J}_{3}\Big[ (1-5c^2)s\cos\theta -\frac{5}{2}s^3\cos3\theta\Big] \\ \nonumber
&& 
+\tilde{J}_{4} \Big[5(1-7c^2)s^2\sin2\theta -\frac{7}{2}s^4\sin4\theta\Big] +\mathcal{O}(e) \Big\}
\\[1ex] \label{d2Z}
\delta_2\Theta &=& \Theta\,\Big\{
-\left[\frac{1}{4}(7-25c^2)+6(1-3c^2)\kappa\right]s^2-\Big[\frac{3}{2}(1-9c^2)+(4-44c^2)\kappa\Big]s^2\cos2\theta \\ \nonumber
&& -\sigma(2-28c^2)s^2\sin2\theta +\frac{3}{4}s^4\cos4\theta-\frac{3}{2}\sigma s^4\sin4\theta \\ \nonumber
&& +\frac{p}{\alpha}\tilde{J}_{3} \Big[
\frac{3}{2}(1-5c^2)s\Big(\sigma\cos\theta+(2+\kappa)\sin\theta\Big)-\frac{5}{4}(4+9\kappa)s^3\sin3\theta+\frac{15}{4}\sigma s^3\cos3\theta\Big] \\ \nonumber
&& -\tilde{J}_{4} \Big[
\frac{5}{2}(1-7c^2)s^2 \Big(2\sigma\sin2\theta+(1+4\kappa)\cos2\theta\Big) -\frac{7}{8}(5+16\kappa)s^4\cos4\theta \\ \nonumber
&& -\frac{7}{2}\sigma s^4\sin4\theta \Big] 
+\mathcal{O}(e^2) \Big\}
\end{eqnarray}

\section*{Acknowledgemnts}

Part of this research has been supported by the Government of Spain (Projects ESP2013-41634-P and ESP2014-57071-R of the Ministry of Economic Affairs and Competitiveness). This research has made use of NASA's Astrophysics Data System Bibliographic Services.

\end{document}